%%%%%%%%%%%%%%%%%%%%%%%%%%%%%%%%%%%%%%%%%%%%%%%%%%%%%%%%%%%
%     GENERALIZED MONOTONE SCHEMES, 
%     EXTREMUM PATHS AND DISCRETE ENTROPY CONDITIONS 
%					
%  Philippe G. LeFloch and Jian-Guo Liu	
%  for publication in Mathematics of Computation
%%%%%%%%%%%%%%%%%%%%%%%%%%%%%%%%%%%%%%%%%%%%%%%%%%%%%%%%%%%
\input amstex
\documentstyle{amsppt}
\magnification 1200
\TagsOnRight
\NoBlackBoxes
%%%%%\hsize=8.10 truein
%%%%%\hcorrection{-.2in}
\vsize=9.0 truein
%%%%%\font\small=cmr8 at 8 pt 
\def\RR{R\!\!\!\!\!I~}
\def\ZZ{\bold Z}  
\def\PP{\Cal P} 
\def\del {{\partial}}
\def\eps{\varepsilon}
\def\minmod{\text{minmod}}
 
\define \refBouchutBourdariasPerthame {1}
\define \refBrenierOsher     {2}
\define \refConwaySmoller    {3} 
\define \refCoquelLeFloch    {4}
\define \refColella          {5}
\define \refCrandallMajda    {6}
\define \refDafermos         {7}
\define \refEngquistOsher    {8}
\define \refFilippov         {9}
\define \refGlimmLax         {10}
\define \refGodlevskiRaviart {11}
\define \refGoodmanLeVeque   {12}
\define \refHartenone        {13}
\define \refHartentwo        {14}
\define \refHartenEngquistOsherChakravarthy {15}
\define \refHartenHymanLax   {16}
\define \refJiangShu         {17}
\define \refKeyfitz          {18}
\define \refKruzkov     {19}
\define \refLaxone      {20}
\define \refLaxtwo      {21}
\define \refLaxWendroff {22}
\define \refvanLeerone  {23}
\define \refvanLeertwo  {24} 
\define \refLeFlochLiu  {25}
\define \refLeroux      {26} 
\define \refLionsSougadinisone  {27}
\define \refLionsSougadinistwo  {28}
\define \refNessyahuTadmor      {29}
\define \refNTT                 {30}
\define \refOsherone  {31}
\define \refOshertwo  {32}
\define \refOsherTadmor{33}
\define \refShu       {34} 
\define \refTadmorone {35}
\define \refTadmortwo {36}
\define \refVolpert   {37}
\define \refYangone   {38}
\define \refYangtwo   {39}
%%%%%%%%%%%%%%%%%%%%%%%%%%%%%%%%%%%%%%%%%%%%%%%%%%%
\topmatter

\

\

\title 
GENERALIZED MONOTONE SCHEMES,
DISCRETE PATHS OF EXTREMA,
AND
DISCRETE ENTROPY CONDITIONS 
\endtitle
\footnote""{\, Published as :  Math. of Comput. 68 (1998), 1025--1055.} 
\author
Philippe  G. LeFloch
\footnote"$^1$"{ \, 
Ecole Polytechnique, Palaiseau, France. CURRENT ADDRESS : 
 Laboratoire Jacques-Louis Lions,
 Centre National de la Recherche Scientifique,
 Universit\'e de Paris 6, 4, Place Jussieu, 75252 Paris, France.
 E-mail address: LeFloch\@ann.jussieu.fr. 
 \newline 
 P.G.L. was 
supported in parts by the Centre National de la Recherche Scientifique, 
and by the National Science Foundation under grants DMS-88-06731, 
DMS 94-01003 and DMS 95-02766, and a Faculty Early Career Development 
award (CAREER). \hskip9.cm  \hfill} 
and Jian-Guo  Liu
\footnote"$^2$"{ \, Temple University, Philadelphia, USA. CURRENT ADDRESS : 
 Department of Mathematics,
University of Maryland,
College Park, MD 20742-4015, USA.  
\newline 
J.G.L. was partially supported by DOE grant DE-FG02 88ER-25053. 
\hskip8.cm 
\hfill}
\endauthor 
\rightheadtext{Generalized monotone schemes}
\leftheadtext{P.G. LeFloch and J.G. Liu}

\endtopmatter
%\abstract 

ABSTRACT. Solutions to conservation laws satisfy the monotonicity property: 
the number of local extrema is a non-increasing function of time, 
and local maximum/minimum values decrease/increase monotonically in time. 
This paper investigates this property from a numerical standpoint. 
We introduce a class of fully discrete in space and time, 
high order accurate, difference schemes, 
called {\it generalized monotone schemes\/}. 
Convergence toward the entropy solution is proven via a new
technique of proof, assuming that the initial data has a 
finite number of extremum values only, and the flux-function is 
strictly convex. 
We define {\it discrete paths of extrema\/} 
by tracking local extremum values in the approximate solution. 
In the course of the analysis we establish the {\it pointwise 
convergence\/} of the trace of the solution along 
a path of extremum. As a corollary, we obtain a proof of 
convergence for a MUSCL-type scheme being second order 
accurate away from sonic points and extrema. 
 
%-----------------------------------------------------------------
\document
\vskip.3cm

\vfill
\eject

%-----------------------------------------------------------------
\heading{1. Introduction}
\endheading 

This paper deals with entropy solutions to the Cauchy problem 
for a one-dimensional scalar conservation law:   
$$
\del_t u + \del_x f(u)  = 0,\qquad 
u(x,t) \in \RR,\,\, t > 0,\,\, x \in \RR,
\tag 1.1
$$
$$
u(x,0)  = u_0(x),\qquad x \in \RR,
\tag 1.2
$$
where the flux $f: \RR \to \RR$ is a given function of class $C^2$ and
the initial data $u_0$ belongs to the space $BV(\RR)$ 
of all integrable functions of bounded total variation. 
For the main result of this paper, we assume that
$$
f \, \text{ is a strictly convex function }
\tag 1.3
$$
and
$$
u_0 \, \text{ has a locally finite number of extrema.}
\tag 1.4
$$
Solutions to conservation laws are generally discontinuous and 
an entropy criterion is necessary to single out a unique solution. 
We refer the reader to Lax \cite{\refLaxone, \refLaxtwo}
for background on nonlinear hyperbolic equations and the entropy criterion. 

As is well-known \cite{\refKruzkov, \refVolpert}, problem (1.1)-(1.2) admits
a unique entropy solution $u$ 
in the space $L^{\infty}\bigl(\RR_+,BV(\RR)\bigr)$ and 
$\text{Lip}\bigl(\RR_+;L^1(\RR)\bigr)$.
(This result holds without the restriction (1.3)-(1.4)
and for multidimensional equations as well.) 
To select the solution one can use the distributional entropy inequality 
$$
\del_t U(u) + \del_x F(u) \le 0,
\tag 1.5
$$
where the (Lipschitz continuous) function 
$(U,F): \RR \to \RR^2$ is a convex entropy-entropy flux pair, i.e.~$U$ is a 
convex function and $F'=U'f'$. 
One can also use the Lax shock admissibility inequality
$$
u(x-,t) \ge u(x+,t)
\tag 1.6
$$
for all $x$ and $t$. 
Since $f$ is convex and $u$ has bounded variation, 
a single entropy $U$ in (1.5) is sufficient to ensure uniqueness, 
and (1.6) is {\it equivalent\/} to (1.5). 
At the discrete level, however, conditions   
(1.5) and (1.6) leads  
two drastically different notions of consistency with the entropy 
criterion for a difference scheme.  
In the present paper we will be using {\it both\/} conditions.   
Indeed in some regions of the $(x,t)$ plane
it is easier to use a discrete version of (1.5),  
while in other regions (1.6) is more adapted. 

We are interested in conservative discretizations of 
problem (1.1)-(1.2) in the sense of Lax and Wendroff \cite{\refLaxWendroff}. 
The monotone schemes and, more generally, the E-schemes are large classes 
of schemes (including the Godunov and Lax-Friedrichs schemes) 
for which a convergence analysis is available. 
The main point is that monotone schemes 
satisfy a discrete version of the entropy inequality (1.5)
(see (2.9) below). However they turn out to be 
first order accurate only, and so have the disadvantage 
to introduce a large amount of numerical viscosity that spreads the  
discontinuities over a large number of computational cells. 

The proof of convergence of the monotone schemes and E-schemes
is based on Helly's and Lax-Wendroff's theorems.  
See \cite{\refConwaySmoller, \refCrandallMajda, \refEngquistOsher, 
\refGodlevskiRaviart, 
\refHartenHymanLax, \refLeroux, \refOsherone, \refTadmorone, 
\refvanLeerone} and the references therein. 
The result holds even for multidimensional equations 
and/or when irregular (non-Cartesian) meshes are used. 

To get high-order accurate approximations, 
it is natural to proceed from analytical properties
satisfied by the entropy solutions to (1.1), formulate them 
at the discrete level, and so design large
classes of high-order difference schemes. One central contribution 
in this direction is due to Harten
\cite{\refHartenone, \refHartentwo}, who introduced the concept 
of TVD schemes, for ``Total Variation Diminishing''. 
Harten shows that conservative, consistent, TVD schemes 
necessarily converge to a weak solution to (1.1). Moreover such 
schemes possess sharp numerical shock profiles with no spurious oscillations. 
However, Harten's notion of TVD scheme is weaker than the 
notion of monotone scheme, and a TVD 
need not converge to the {\it entropy\/} solution.  
The aim of this paper is precisely to single out a subclass of TVD schemes, 
refining Harten's notion, 
that are both high-order order accurate and entropy satisfying, 
cf.~Definition 2.2 below.

A very large literature is available on the actual design of 
second-order shock-capturing schemes. One approach 
to upgrading a first order scheme was proposed by van Leer 
\cite{\refvanLeerone, \refvanLeertwo}: the MUSCL scheme 
extends the Godunov scheme by replacing the piecewise 
constant approximation in the latter with a  
piecewise affine approximation. The heart of the matter is to 
avoid the formation of spurious oscillations near discontinuities. 
This is achieved by van Leer via the so-called min-mod limitor.

Other classes of high order schemes 
have been built from the maximum principle 
and a monotony condition: see the classes of ENO and UNO schemes 
proposed by Harten, Engquist, Osher, and Chakravarthy
\cite{\refHartenEngquistOsherChakravarthy}. 
In this paper we concentrate attention on the MUSCL scheme, but 
our main convergence theorem, Theorem 2.3 below, should also 
be applicable to other schemes.

After Osher's pioneering result \cite{\refOshertwo} and the systematic study  
of high-order schemes by Osher and Tadmor \cite{\refOsherTadmor} 
(which however required a technical condition on the slopes of the affine 
reconstructions), Lions and Sougadinis \cite{\refLionsSougadinisone, 
\refLionsSougadinistwo} and, independently, 
Yang \cite{\refYangone} established the convergence of the MUSCL scheme. 
Both proofs apply to an arbitrary BV initial data (so (1.4) is not assumed); 
the flux-function is assumed to be convex in \cite{\refYangone} 
and strictly convex in \cite{\refLionsSougadinistwo}. 
In both papers significantly new techniques of proof are introduced
by the authors. 
In \cite{\refYangone}, Yang develops a method for tracking local 
extremum values.  
In \cite{\refLionsSougadinistwo}, Lions and Sougadinis elegantly  
re-formulate the MUSCL scheme   
at the level of the Hamilton-Jacobi equation associated with (1.1) 
and rely on Crandall-Lions' theory of viscosity solutions for such equations. 
A large class of difference approximations
is treated in \cite{\refLionsSougadinistwo}. 
Techniques in both papers are restricted in an essential way 
to semi-discrete schemes, 
in which the time variable is kept continuous. 
In \cite{\refLeFlochLiu} the present authors announced a proof of 
the convergence of a class of fully discrete schemes that include 
the MUSCL scheme. Independently, Yang \cite{\refYangtwo} also extended 
his approach to a large class of fully discrete methods. 

Several other works deal with the convergence of 
the van Leer's scheme or variants of it. A discrete version of inequality 
(1.6) was established by Brenier and Osher \cite{\refBrenierOsher} and  
Goodman and LeVeque \cite{\refGoodmanLeVeque}, the latter dealing with 
both first and second order methods. 
Nessyahu, Tadmor, and Tamir \cite{\refNTT} establish both the convergence 
and error estimates for a variety of Godunov-type schemes. 
Various Approaches to upgrading Lax-Friedrichs 
scheme are actively developed by Tadmor and co-authors; see, for 
instance, Nessyahu and Tadmor \cite{\refNessyahuTadmor}. 
An extensive discussion of the discretization of the entropy inequality 
(1.5) is found in Osher and Tadmor \cite{\refOsherTadmor} 
and the many references cited theirein. See also 
Bouchut, Bourdarias and Perthame \cite{\refBouchutBourdariasPerthame},
and Coquel and LeFloch \cite{\refCoquelLeFloch}.

The objective of the present paper is to provide a framework to 
prove the convergence of high order accurate and
fully discrete difference schemes. We will strongly rely on 
a property shared by all entropy solutions to (1.1):
the {\it monotonicity property\/}. 
Given an arbitrary entropy solution $u$, the number of extrema in $u(t)$ is 
a non-increasing function of $t$, and local maximum/minimum values 
decrease/increase monotonically in time. 
(See the Appendix for rigorous statements.)  
This property was studied first by Harten \cite{\refHartenone}
from the numerical standpoint, in order to arrive at his notion 
of TVD scheme. It was also essential 
in \cite{\refLionsSougadinistwo, \refYangone} and in 
Tadmor \cite{\refTadmortwo}.

Observe that monotone or TVD schemes 
do not necessarily satisfy the monotonicity 
property. For instance the Lax-Friedrichs scheme may increase the 
number of extremum values! 
This motivates the introduction of a subclass of TVD schemes, 
guaranteeing this property together with the high order of accuracy. 
The scheme then closely mimics an important property of 
the solutions to the continuous equation.  Further requirement 
is necessary to ensure that the scheme is entropy-satisfying. 

Based on the monotonicity property, we thus introduce the notion of 
{\it generalized monotone scheme\/} 
(Definition 2.2) characterized by three conditions: 
maximum principle, entropy consistency, and local behavior at extrema. 
The (first order) Godunov method and the (second order) 
MUSCL method are prototype examples. 

Our first requirement is 
a (strong) version of the {\it local maximum principle\/}, 
which will allow us to show a discrete analogue of the monotonicity property. 
In particular our condition prevents the formation of spurious numerical 
oscillations and the scheme is TVD in the sense of Harten. 

The second requirement is motivated by the following observation made by Osher 
for semi-discrete schemes \cite{\refOsherone}:  
in the non-decreasing parts of the approximate solutions, 
it is possible to simultaneously achieve second order accuracy and 
the existence of one cell entropy inequality 
--i.e.~a discrete analogue of (1.5).
Our second condition therefore requires 
a {\it cell entropy inequality\/} in the {\it non-decreasing\/} regions.

Finally, in order to prevent cusp-like behavior near extrema, which may lead to 
entropy violating discontinuities, the scheme should be well-behaved 
near local extrema. To this end, 
we introduce a condition referred to as 
the {\it quadratic decay property at local extrema\/}. 
It reflects the nonlinear behavior of the numerical flux-function near
extrema, and assumption (1.3) again is essential. For simplicity, 
it is also assumed that the scheme reduces to a three points, first order 
scheme at extrema.

The rest of this section is devoted to comments upon the proof of 
our main result that any generalized monotone scheme 
converges to the entropy solution of (1.1)-(1.2) provided (1.3)-(1.4) holds. 
Our approach was driven by Yang's paper \cite{\refYangone} on semi-discrete 
schemes. However, our technical arguments differ 
substantially from the ones in \cite{\refYangone}. 
In particular we restrict attention 
to initial data having a locally finite number 
of extrema, making the tracking of paths of extremum values 
almost a trivial matter. The main part of our proof is studying the 
convergence of the traces of the approximate solution
along it. We make use of the quadratic decay property 
above to exclude the formation of ``cusp'' near extremum values, 
which could lead to entropy-violating shock. This is the main 
contribution of this present paper. 
We do not believe that the extension of our proof to arbitrary BV data is
straightforward. The result we obtain seems satisfactory however since 
condition (1.4) covers ``generic'' initial data.

The proof distinguishes between the non-increasing parts and
non-decreasing parts of the solution, 
and is based on several observations as follows. 
We use the notation $u^h$ for the approximate solutions,
$v$ for the limiting solution, and $h$ for the mesh size.  

First of all, the strong maximum principle ensures that the $u^h$'s 
are total variation diminishing in time, so of uniformly bounded
total variation. The strong maximum
principle enables to easily define a discrete path and each time
step the discrete path move at most one grid points. As a consequence, 
the discrete path is Lipschitz continuous. (This is a major
difference between the present paper and Yang's paper in which the 
construction of the path and the limiting paths are much more involved.)

By Helly's Theorem, the scheme 
converges in the strong $L^1$ topology
to a limiting function, say $v$, which according 
to Lax-Wendroff's theorem is a weak solution to (1.1)-(1.2).
It remains to prove that $v$ is {\it the\/} entropy solution.
We show that the strong maximum principle in fact implies a 
discrete analogue of the monotonicity property.
Relying on (1.4), we construct a (locally) finite family of 
Lipschitz continuous paths in the plane
by tracking the local extrema in $u^h$. 
The paths are shown to converge in the uniform topology to limiting curves.

Next we make the following two observations.
On one hand, in a non-increasing region for $u^h$, the function $v$ is also 
non-increasing, and so can only admit non-increasing jumps. Thus $v$ 
satisfies the Lax shock admissibility inequality (1.6) in the 
non-increasing regions. On the other hand, 
a discrete cell entropy inequality, by assumption, holds 
in the non-decreasing regions of $u^h$. 
So $v$ satisfies (1.5) in the non-decreasing regions. 

It remains to prove that $v$ has only {\it non-increasing jumps\/}
 along any path
of extrema. This is the most interesting part of the proof.
Let $\psi^h$ be an approximate path of extrema, and let $\psi$ be its uniform limit. 
The path $\psi$ is the boundary separating two regions
where the analysis in the paragraph above applies. Indeed 
(1.5) holds in the side where $v$ is non-increasing
and (1.6) holds in the side where $v$ is non-decreasing. 
A specific proof must be provided to determine the behavior of $v$ 
{\it along the path\/}. We analyze 
the entropy production in a small region of the plane limited   
on one side by the path $\psi^h$. 
In the course of this proof we derive a uniform bound for the 
time integral of the local oscillation in space 
along the path $\psi^h$, which 
is a direct consequence of the quadratic decay property mentioned above. 
For the sake of simplicity, we use here the assumption that the scheme 
reduces to a three points, first order scheme at extrema.

Note that assumption (1.3) is not used in the construction of the 
extremum paths, but is essential in the
convergence analysis which strongly relies upon (1.6), only valid for
convex fluxes. 

Our analysis shows that, for a class of difference schemes, certain approximate 
generalized characteristics --those issued from an extremum point of the initial data-- 
can be constructed for scalar conservation laws with convex flux. 
Constructing approximate generalized characteristics issued from an
arbitrary point remains a challenging open problem. Recall that,
for the random choice scheme, 
Glimm-Lax \cite{\refGlimmLax} did obtain a general theory of approximate
generalized characteristics (applicable to systems, as well).

The organization of this paper is as follows. 
In Section 2, we define the class of generalized monotone 
schemes and state the main result
of convergence, cf.~Theorem 2.3.  
Section 3 contains the proof of the main result.
In Section 4, we apply Theorem 2.3 to the MUSCL scheme and provide 
some additional remarks.

%--------------------------------------------------------------
\heading{2. Generalized Monotone Schemes}
\endheading 

This section introduces a class of TVD schemes which 
are built to closely mimic an essential property of the 
entropy solutions to (1.1), i.e.~the monotonicity property: 
the number of local extrema is a non-increasing function of time, 
and local maximum/minimum values decrease/increase monotonically in time. 
See Appendix for a precise statement. 
We investigate here this property from a numerical standpoint. 
Monotone and TVD schemes actually do not necessarily satisfy this 
property (see Section 4 for an example), and a more restricted 
class, the {\it generalized monotone schemes\/}, is natural.

We consider a $(2k+1)$--point difference scheme in conservation form 
for the approximation of (1.1)-(1.2): 
$$
u^{n+1}_j = u^n_j - \lambda \bigl( g^n_{j+{1/2}}-g^n_{j-{1/2}}\bigr), 
\quad n \ge 0, \, j \in \ZZ, 
\tag 2.1
$$
where we use the notation
$g^n_{j+{1/2}} = g\bigl(u^n_{j-k+1}, \dots, u^n_{j+k}\bigr)$
and $\lambda=\tau/h$ for the
ratio of the time-increment $\tau$ by the space-increment $h$. We set
$t_n=n\tau$, $x_j = jh$, and $x_{j+{1/2}}=(j+1/2)h$.  The
value $u^n_j$ presumably is an approximation of the exact solution at the point
$(x_j, t_n)$. As is usual, the numerical flux $g:
\RR^{2k}\to\RR$ is assumed to be locally Lipschitz continuous and
consistent with $f$, i.e.~$g(v,\dots,v)=f(v)$ for all $v$. 
Note that $g$ may depend on $\lambda$. 
For definiteness, we set
$$
u^0_j=\frac1h\int_{x_{j-{1/2}}}^{x_{j+{1/2}}} u_0(y)\ dy.
\tag 2.2
$$
This is sufficient for second order accuracy. 
For higher orders, one
should use a Runge-Kutta time-step method (Shu \cite{\refShu}). 
We also define the piecewise constant function $u^h: \RR\times\RR_+\to\RR$ by 
$$
u^h(x,t) = u^n_j, \quad t_n\le t <t_{n+1}, \, \, x_{j-{1/2}} \le x < x_{j+{1/2}}. 
\tag 2.3
$$
By construction $u^h$ is a right continuous function. 
For simplicity we assume the following CFL restriction:   
$$
\lambda \, \sup_{v} |f'(v)| \le 1/4.
\tag 2.4
$$
Several of the properties below would still hold if, in (2.4), one
replaces $1/4$ by 1 although the proofs then become less clear geometrically. 

The (first order) Godunov scheme is based on exact solutions to (1.1), 
and is a good prototypical example to lead us to defining a class of 
(high order) schemes consistent with the monotonicity property.

A main ingredient is the Riemann
problem. Given two states $v$ and $w$, we define 
$R(.;v,w)$ to be the entropy solution to  
(1.1)-(1.2) with, here,  
$$
u_0(x) = v \quad \text { if } x<0,   \qquad w \quad \text { if } x>0. 
$$
As is well-known, $R(\cdot;v,w)$ depends on the self similarity
variable $x/t$ only, and is given by a closed formula. If $v\le w$, $R$ is a
rarefaction wave and, if $v>w$, a shock wave. More important 
$$
R(\frac x t;v,w) \text{ is a monotone function connecting } 
v \text{ to } w.
\tag 2.5
$$

In Godunov scheme, 
one solves Riemann problems and, at each time level, one projects the
solution on the space of piecewise constant functions. If  $u^h(t_n)$ is
known, let $\tilde u(x,t)$ for $t\ge t_n$ be the entropy solution to
(1.1) assuming the Cauchy data
$$
\tilde u(t_n+) = u^h(t_n+).
$$
Since $u^h(t_n+)$ is a piecewise constant function, $\tilde u$ is obtained 
explicitly by glueing together the Riemann solutions $R(.;u^n_j,u^n_{j+1})$. 
In view of (2.4), there is no interaction
between two nearby solutions, at least for $t \in (t_n,t_{n+1})$. 
Set    
$$
u^h(x,t_{n+1}+) = \frac 1 h \int_{x_{j-1/2}}^{x_{j+1/2}}
\tilde u(y,t_{n+1}) \, dy, \quad x_{j-{1/2}} \le x < x_{j+{1/2}}.
$$
Using the conservative form of (1.1), the scheme can be
written in the form (2.1) with $k=1$ and $g=g_G$ given by 
$$
g_G(v,w) = f(R(0+;v,w)) \quad \text{ for all } v \text{ and } w.
\tag 2.6
$$

The proofs of (2.7)--(2.10) stated below are classical matter; 
e.g.~\cite{\refConwaySmoller, \refCrandallMajda, \refHartenHymanLax, 
\refLeroux}. 
On one hand, one can think of the Godunov scheme 
geometrically as a two-step method:
a marching step based on exact (Riemann) solutions
and an $L^2$ projection step. 
In view of (2.5) and the monotonicity property 
of the $L^2$ projection, one easily have 
a simple geometrical proof of the properties listed below. 
On the other hand, an algebraic approach is 
based on the explicit formula deduced from (2.6).

The Godunov scheme is monotone: the function $g_G$ is 
non-decreasing with respect to its first argument,
and non-increasing with respect to its second one. 
This property implies that the scheme is 
{\it monotonicity preserving\/}, i.e., 
$$
\aligned
& \text{ if } u_{j_1}^n,  u_{j_1+1}^n, \dots, u_{j_2}^n 
\text{ is a non-increasing (resp. non-decreasing) sequence }\\
& \text{ for some indices } j_1< j_2, 
\text{ so is the sequence } 
u_{j_1+1}^{n+1},  u_{j_1+2}^{n+1}, \dots, u_{j_2-1}^{n+1}.\cr
\endaligned
\tag 2.7
$$
The Godunov scheme satisfies the {\it local maximum principle\/}, 
 i.e., 
$$
\min\big(u^n_{j-1}, u^n_j, u^n_{j+1}\big) 
\le u^{n+1}_j \le \max \big(u^n_{j-1}, u^n_j, u^n_{j+1}\big) 
\tag 2.8
$$
for all $n \ge 0$ and $j \in \ZZ$. 
In fact (2.7) and (2.8) are shared by both steps in the Godunov scheme.
It will convenient to us to rewrite (2.8) in term of the 
jumps of $u^h$ at the endpoints of a cell:      
$$
\min\big(u^n_{j+1} - u^n_j, u^n_{j-1} - u^n_j,0\big) 
\le u^{n+1}_j -  u^n_j 
\le \max\big(u^n_{j+1} - u^n_j, u^n_{j-1} - u^n_j,0\big).
\tag 2.8'
$$

Any monotone scheme --in particular the Godunov scheme--
satisfies a discrete analogue of (1.5)
for every convex entropy pair $(U,F)$: 
$$
U(u^{n+1}_j) -  U(u^n_j)  - \lambda 
\bigl( G^n_{j+{1/2}} - G^n_{j-{1/2}} \bigr) \le 0, \quad n\ge 0, \, j \in \ZZ. 
\tag 2.9
$$
In (2.9), 
$ G^n_{j+{1/2}} = G(u^n_{j-k+1}, \dots, u^n_{j+k})$, and 
$G$ is a numerical entropy flux consistent with $F$, that is 
$G(v,v\dots,v) = F(v)$ for all $v$.

Finally concerning the local behavior of $u^h$ in the neighborhood
of local extrema, it is known that, say for local maximum,  
$$
u_j^{n+1} \le u_j^n. 
\tag 2.10 
$$
A similar property holds for local minima. 

In fact the classical properties (2.8) and (2.10)
can be improved as follows.

\proclaim{Proposition 2.1} 
Under assumptions $(1.3)$ and $(2.4)$, 
the Godunov scheme satisfies the following two properties: 
\roster
\item the {\rm strong local maximum principle}: 
$$
\frac 1 2 \min\big(u^n_{j+1} - u^n_j, u^n_{j-1} - u^n_j,0\big) 
\le u^{n+1}_j -  u^n_j 
\le \frac 1 2 \max\big(u^n_{j+1} - u^n_j, u^n_{j-1} - u^n_j,0\big)
\tag 2.11
$$
\item and the {\rm quadratic decay property at local extrema}, 
that is e.g.~for a maxima: 
$$
\text{If } \, u_j^n \, \text{ is a local maximum value }, \quad 
u_j^{n+1} \le 
u^n_j - \alpha \min_\pm\bigl((u^n_j -u^n_{j \pm 1})^2\bigr). 
\tag 2.12
$$
\endroster
with $\alpha= \lambda \inf f'' /2$.  $\hfill \square$
\endproclaim

The proof of (2.11) is straightforward from a geometrical standpoint. 
If also follows from Proposition 4.1 established later in 
Section 4. 
Note that the coefficient $1/4$ in (2.4) is essential for (2.11) to hold.  
Observe that (2.11) is stronger than (2.8)-(2.8')
and controls the time-increment (i.e., $u^{n+1}_j- u^n_j$)
of the solution in the cell $j$ in term of the jumps at the endpoints: 
the values $u^n_j$ evolves ``slowly'' as $t_n$ increases.  
As we shall see, this property implies that the scheme satisfies a discrete
analogue of the monotonicity property. 

Estimate (2.12) is stronger than (2.10)
and shows that the decrease/increase of a maximum/minimum is 
controlled by the quadratic oscillation of $u^h$ nearby this extremum. 
It is a truly nonlinear property of the Godunov flux.  
It will be used below to prove that cusp 
can not form near extremum points. For convenience, 
the proof of (2.12) is postponed to Section 4, 
where second-order approximations are treated as well.   

In \cite{\refGoodmanLeVeque}, Goodman and LeVeque derive for the Godunov 
method a discrete version of the Oleinik entropy inequality. 
In particular, this shows that the Godunov scheme spreads rarefaction 
waves at the correct rate. Our estimate (2.12) is, at least in spririt,  
similar to this spreading estimate, and expresses the spreading of 
extremum values. 

We are now ready to introduce a class of high-order schemes 
based on the properties derived in Proposition 2.1.

\proclaim{Definition 2.2} 
The scheme $(2.1)$ is said to be a {\rm generalized monotone scheme}  
if any sequence
$\big\{u^n_j\big\}$ generated by $(2.1)$ satisfies the 
following three conditions:    
\roster
\item the strong local maximum principle $(2.11)$,
\item the cell entropy inequality $(2.9)$ for one strictly 
convex pair $(U,F)$ in any non-decreasing region, 
including local extrema,   
\item the quadratic decay property at local extrema $(2.12)$ 
for some constant $\alpha >0$. 
\endroster
It is also assumed that the numerical flux and the numerical entropy flux 
are essentially two-point functions at local extrema. $\hfill \square$
\endproclaim

According Proposition 2.1, the (first order) Godunov scheme belongs to the class
described in Definition 2.2.
Section 4 will show that there exist high order accurate schemes
satisfying the conditions in Definition 2.2. 
Our main convergence result is: 

\proclaim{Theorem 2.3} 
Let $(2.1)$ be a generalized monotone scheme. Assume that assumptions $(1.3)$-$(1.4)$
hold together with $(2.4)$. Then the scheme $(2.1)$  
\roster
\item is $L^\infty$ stable, i.e., 
$$
\inf_{l \in \ZZ} u_l^n \le u_j^{n+1}  \le \sup_{l \in \ZZ} u_l^n, 
\qquad  n \ge 0, \, j \in \ZZ, 
\tag 2.13
$$
\item is total variation diminishing, i.e., 
$$
\sum_{j \in \ZZ} \bigl| u^{n+1}_{j+1} - u^{n+1}_j \bigr| 
\le \sum_{j \in \ZZ}  \bigl| u^n_{j+1} - u^n_j \bigr|, \qquad n \ge 0, 
\tag 2.14
$$
\item and converges in the $L^p_{\text{loc}}$ strong topology for all 
$p \in [1,\infty)$ to the entropy solution of $(1.1)$-$(1.2)$. $\hfill \square$
\endroster 
\endproclaim

The proof of Theorem 2.3 is given in Section 3. 

Theorem 2.3 is satisfactory for a practical standpoint.  
Suppose that $u_0$ is an arbitrary BV function, 
and we wish to compute an approximation to the
solution $u$ of (1.1)-(1.2) of order $\eps>0$ in the $L^1$ norm. 
Let us determine first an approximation of $u_0$, say $u_{0,\eps}$, 
that has a finite number of local extrema and such that      
$$
\| u_{0,\eps} - u_0 \|_{L^1(\RR)} \le \eps.
$$
Applying a generalized monotone scheme to the initial condition 
$u_{0,\eps}$ yields an approximate solution $u^h_\eps$ that,
in view of Theorem 2.3, satisfies   
$$
\|u^h_\eps - u_\eps \|_{L^1(\RR)} \le o(h) \le \eps
$$
for $h$ is small enough, 
where $u_\eps$ is the entropy solution associated with the initial
condition $u_{0,\eps}$. Since the semi-group of solutions associated with 
 (1.1) 
satisfies the $L^1$ contraction property, one has 
$$
\|u_\eps - u \|_{L^1(\RR)} \le   \| u_{0,\eps} - u_0 \|_{L^1(\RR)} \le \eps, 
$$
and therefore
$$
\align
\|u^h_\eps - u\|_{L^1(\RR)} 
& \le  \| u^h_\eps - u_\eps \|_{L^1(\RR)} + \|u_\eps -u \|_{L^1(\RR)}  
\\
& \le 2 \, \eps.
\endalign
$$

%--------------------------------------------------------------------
\heading{3. Convergence Analysis}
\endheading 

The proof of Theorem 2.3 is decomposed into several lemmas, 
Lemmas 3.1--3.12. For the whole of this section,
we assume that the hypotheses made in Theorem 2.3 are
satisfied. 

We introduce first some notation and terminology. 
We call $u_j^n$ a local maximum or a local minimum if there 
exist two indices $j_*$ and $j^*$ with $j_* \le j \le j^*$ such that 
$$
u_{j_*}^n = u_{j_*+1}^n = \dots =  u_{j^*}^n 
> \max (u_{j_*-1}^n, u_{j^*+1}^n)
$$
or
$$
u_{j_*}^n = u_{j_*+1}^n = \dots =  u_{j^*}^n 
< \min (u_{j_*-1}^n, u_{j^*}^n).
$$
In such a case, there is no need to distinguish between  
the extrema $u_{j_*}^n$, $u_{j_*+1}^n$,
$\dots$, $u_{j^*}^n$. 
Based on the strong maximum principle (2.11), we show in
Lemmas 3.1 and 3.2 that the scheme satisfies
a discrete form of the monotonicity property. 
We construct a family of paths in the $(x,t)$-plane
by tracing in time the points where the approximate solution $u^h(t)$ 
achieves its local extremum values. 
One difficulty is proving that the interaction of two (or more)
paths does not create new paths, so the total number of paths at any given time 
remains less or equal to the initial number of local extrema in $u_0$. 
In passing we observe 
that an extremum point moves one grid point at each time-step, at most.

\proclaim{Lemma 3.1} 
For some $j_* < j^*$, suppose that the sequences 
$\,u^n_{j_*-3}$, $u^n_{j_*-2}$, $u^n_{j_*-1}$, $u^n_{j_*}$ 
and 
$\,u^n_{j^*}$, $u^n_{j^*+1}$, $u^n_{j^*+2}$, $u^n_{j^*+3}$ 
are two monotone sequences, no
specific assumption being made on the values $u^n_j$,
$j_*\le j \le j^*$.
Then the number $\nu'$ of extrema in the sequence 
$$
S_{n+1} := \big(u^{n+1}_j\big)_{j_*-2 \le j \le j^*+2} 
$$
is less or equal to the number $\nu$ of extrema in
$$
S_n := \big(u^n_j\big)_{j_*-2 \le j \le j^*+2}. 
$$
When $\nu'\ge 1$, 
there exists a one-to-one correspondence between $\nu'$
local extrema of $S_n$  
and the $\nu'$ local extrema of $S_{n+1}$ with the following property:
if a maximum/minimum $u^n_j$ is associated 
with a maximum/minimum $u^{n+1}_{j'}$, then 
$$
|j' - j| \le 1 \quad \text{ and } \quad u^{n+1}_{j'} \le u^n_j, 
\text{ resp. } \, u^{n+1}_{j'} \ge u^n_j.
\tag 3.1 
$$ 
$\hfill \square$
\endproclaim

\proclaim{Proof} \rm 
We distinguish between various cases
depending on the number of local extrema in the sequence 
$S_n$ and construct the one-to-one correspondence. 

If $S_n$ has no local extremum, for instance is 
non-decreasing, then $S_{n+1}$ 
is also non-decreasing. This indeed follows from inequalities (2.11)
which reduce in this case to
$$
\cdots \leq \frac 1 2 \big(u^n_{j-1} + u^n_j\big) 
\leq u^{n+1}_j  
\leq \frac 1 2 \big(u^n_{j+1} + u^n_j\big) \leq u^{n+1}_{j+1} \leq \cdots 
$$

Consider next the case that $S_n$ has exactly one local extremum, 
say a local maximum at some $u^n_l$. The same argument as above shows that 
the sequences $\big\{u_j^{n+1}\big\}_{j_*\le j\le l-1}$ and 
$\big\{u_j^{n+1}\big\}_{l+1\le j\le j^*}$ are non-decreasing
and non-increasing respectively. Therefore we only need to exclude the case that 
both $u_{l-1}^{n+1}>u_l^{n+1}$ and $u_{l+1}^{n+1}> u_l^{n+1}$, 
which would violate the monotonicity property since  $S_{n+1}$ 
in this case 
would have two local maximum and one local minimum, so two new extrema. 
Indeed assume that the latter would hold, then using (2.11) at the points
$l-1$, $l$, and $l+1$ gives us
$$
u^{n+1}_{l-1} \le \frac 1 2 \big(u^n_{l-1} + u^n_l\big), 
$$
$$
u_l^n  + \frac 1 2 \min\big(u^n_{l\pm 1} - u^n_l\big) 
\le u^{n+1}_l,
$$
and
$$
u^{n+1}_{l+1} \le \frac 1 2 \big(u^n_{l+1} + u^n_l\big), 
$$
which are incompatible with the inequalities $u_{l-1}^{n+1}>u_l^{n+1}$ 
and $u_{l+1}^{n+1}> u_l^{n+1}$. 

Consider now the case that  $S_{n}$ has two local extrema, say one
local maximum at $l$ and one local minimum at $m$ with $l <m$. 
We distinguish between three cases:

If $l<m-2$, then the two extrema can not ``interact'' and 
the arguments before show that the solution at time $t_{n+1}$
has the same properties.
  
If $l=m-2$, the two extrema can interact. Using (2.11) 
at each point $j=l-1, \cdots, l+3$, one gets 
$$
u^{n+1}_{l-1} \le \frac 1 2 \big(u^n_{l-1} + u^n_l\big), 
$$
$$
u_l^n  + \frac 1 2 \min\big(u^n_{l\pm 1} - u^n_l\big) 
\le u^{n+1}_l,
$$
and
$$
u^{n+1}_{l+1} \le \frac 1 2 \big(u^n_{l+1} + u^n_l\big), 
$$
and 
$$
u^{n+1}_{l+1} \le \frac 1 2 \big(u^n_{l+1} + u^n_{l+2}\big), 
$$
$$
u_{l+2}^n  + \frac 1 2 \min\big(u^n_{l+2 \pm 1} - u^n_{l+2}\big) 
\le u^{n+1}_{l+2},
$$
and
$$
u^{n+1}_{l+3} \le \frac 1 2 \big(u^n_{l+3} + u^n_{l+2}\big). 
$$
It is not hard to see that these inequalities imply that  $S_{n+1}$ 
\roster 
\item either has one maximum at either $j=l-1, l, l+1$
and a minimum at $j=l+1, l+2, l+3$, 
\item or is non-decreasing. 
\endroster 
In the first case, we achieve the property we wanted. 
In the second case, there is no extremum at the time $t_{n+1}$. 

This analysis can be extended to the case that several 
extrema can ``interact''; we omit the details. 
Property (3.1) is a consequence of condition (2.4): 
an extremum point can only move up to one grid point at each time step. 
$\hfill \square$
\endproclaim

Consider the initial condition $u_0$ and its approximation $u^h(0)$
defined by $L^2$ projection, cf.~(2.2). 
Locate the minimum and maximum values in the initial data $u_0$. 
For $h$ much smaller than the minimal distance between two 
extrema, $u^h(0)$ has the same number of extrema
as $u_0$ and the same increasing/decreasing
behavior as $u_0$. Indeed, 
there exist indices $J_q^h(0)$ for $q$ in 
a set of consecutive integers $E(u_0)$ depending on $u_0$ but not on $h$, 
such that
$$
\aligned 
&  u^0_j \text{ is non-decreasing for } \, J^h_{2p}(0) \le j \le J^h_{2p+1}(0) \\
&  u^0_j \text{ is non-increasing for } \, J^h_{2p-1}(0) \le j \le J^h_{2p}(0). 
\endaligned 
\tag 3.2 
$$
Those indices are not uniquely determined in the case that $u_0$ is 
constant on an interval associated with a local extremum. 
Since $u_0$ has a locally finite number of local extrema, 
there exists a partition of $\ZZ$ into intervals $(j_*,j^*)$ 
in which the hypothesis of Lemma 3.1 holds. 
It is an easy matter to use the one-to-one correspondence
in Lemma 3.1 and trace forward in time up to time $t_1= \tau$ 
the locations of the extrema in  $(j_*,j^*)$ . 
At each time level a (possibly new) partition of $\ZZ$
is considered and Lemma 3.1 is  used again.
Indeed the values $J_q^h(n+1)$ in Lemma 3.2 below are defined 
from the $J_q^h(n)$'s according to the one-to-one 
correspondence established in Lemma 3.1.
Finally piecewise affine and continuous paths are obtained by connecting 
together the points of local extrema. 
It may happen that the number of extrema 
decreases from time $t_n$ to $t_{n+1}$.  
In such a case, one path, at least, can no longer be further extended in time  
and so, for that purpose,
we introduce a ``stopping time'', denoted by $T_q^h=t_n$. 

The following lemma is established.

\proclaim{Lemma 3.2} 
There exist continuous and piecewise affine curves
$\psi_q^h: [0,T_q^h] \to \RR$ for $q \in E(u_0)$,
passing through the mesh points $\big(x_{J_q^h(n)},t_n\big)$ and 
having the following properties: 
$$
\psi_q^h(t) = x_{J_q^h(n)} + \frac{t-t_n}{\tau} 
\bigl( x_{J_q^h(n+1)}- x_{J_q^h(n)} \bigr), \quad t \in [t_n,t_{n+1}], 
\tag 3.3
$$
for each $n=0,1,2,\dots, N_q^h$ with $T_q^h = N_q^h \tau \le \infty$, 
$$
\psi_q^h \le \psi_{q+1}^h, \qquad |x_{J_q^h(n)} -  x_{J_q^h(n+1)}| \le h; 
\tag 3.4
$$
$$
\aligned
& \text{ there is only a finite number (uniformly bounded w.r.t. h) } 
\text{ of curves  }\psi_q^h \cr
& \text{ on each compact set} \cr
\endaligned
\tag 3.5
$$
and 
$$
\aligned
&  x \in (\psi_{2p}^h(t_n), \psi_{2p+1}^h(t_n)) \mapsto 
   u^h(x,t_n) \,\, \text{ is non-decreasing,} \\
&  x \in (\psi_{2p-1}^h(t_n), \psi_{2p}^h(t_n)) \mapsto 
   u^h(x,t_n) \,\, \text{ is non-increasing.}
\endaligned
\tag 3.6
$$
Furthermore, the functions $w_q^h: [0,T_q^h] \to \RR$ defined by
$$
w_q^h(t) = u_{J_q^h(n)}^n \qquad \text{ for } \qquad t_n \le t < t_{n+1}
\tag 3.7
$$
are non-decreasing if $q$ is even, and non-increasing if $q$ is odd. $\hfill \square$
\endproclaim

\proclaim{Remark 3.3} \rm 1) \ 
The definition (3.3) is not essential.
All the results below still hold if $\psi_q^h$ is replaced 
by any (uniformly) Lipschitz continuous curve 
passing through the mesh points $\big(x_{J_q^h(n)},t_n\big)$. 
As a matter of fact, it is an open problem to show
the strong convergence of the derivatives of approximate paths. 
By comparison, 
for the approximate solutions built by the random choice scheme, 
Glimm and Lax \cite{\refGlimmLax} prove the a.e. convergence
of the first order derivatives of the paths. \par 
\noindent 2) \ Introducing the stopping times $T_q^h$ is necessary. 
At those times, certain paths cross each other 
and their extension in time is not well-defined. 
For instance a path of maximum 
and a path of minimum  can cross and ``cancel out''.
The case of exact solutions (Cf.~the appendix) is simpler in this respect: 
the paths can be defined to be characteristic curves for all times, 
even when they are no longer paths of extrema. \par 
\noindent 3) \ It is not interesting to trace the minimal
(or maximal) paths of extrema in the approximate solution. 
Such paths would not converge to the 
paths obtained in the continuous case.  $\hfill\square$
\endproclaim

By construction, cf.~(3.4), a path ``jumps'' up to 
one grid point at each time-step.
So the slope of a path remains uniformly bounded by $1/\lambda$ and 
the curves $\psi_q^h$ are bounded in the $W^{1,\infty}_{\text loc}$ norm,
uniformly with respect to $h$ and $q$. 
On the other hand, Lemma 3.3 implies that the scheme is TVD so
$TV(u^h(t_n))$ is uniformly bounded.
We thus conclude that the approximate solutions
and the paths of extrema are strongly convergent, as stated in 
the following Lemmas 3.4 and 3.5. 
For simplicity, we keep the same notation for a sequence and 
a subsequence.

\proclaim{Lemma 3.4}  
There exist times $T_q \in [0, +\infty]$ and Lipschitz continuous curves 
$\psi_q:(0,T_q) \to \RR$ such that
$$
T_q^h\to T_q  \quad \text{ as } \, h \to 0,
\tag 3.8
$$
and 
$$
\psi_q^h \to \psi_q \quad 
\text{ uniformly on each compact subset of } C^0([0,T_q)).
\tag 3.9
$$
$\hfill \square$
\endproclaim

\proclaim{Lemma 3.5}  
The sequence $u^h$ satisfies estimates $(2.13)$-$(2.14)$, 
and so is uniformly stable in the $L^\infty([0,\infty), BV(\RR))$ and
${\text{Lip}}([0,\infty), L^1(\RR))$ norms. There exists
a function $v$ in the same spaces such that
$$
u^h(x,t) \to v(x,t)  \qquad\text{ for all times } \,  t\ge 0
\, \text{ and almost every } \, x \in \RR, 
\tag 3.10
$$
and there exist functions $w_q$ in $BV((0, T_q),\RR)$ such that
$$
w_q^h \to w_q \qquad \text{ almost everywhere on } (0, T_q)
\tag 3.11
$$
for all $q \in E(u_0)$. 
$\hfill \square$
\endproclaim

The convergence results (3.10) and (3.11) hold in particular
at each point of continuity of $v(t)$ 
and $w_q$, respectively. 
Introduce now the following three sets, 
which provide us with a partition of the $(x,t)$-plane 
into increasing/decreasing regions for $v$: 
$$
\aligned
& 	{\Omega}_1(v) = \big\{(x,t) \, / \, \psi_{2p_1}(t) < x < \psi_{2p_2+1}(t),  
\, \, p_1 \le p_2, \, t<T_{2p_1}, \, t<T_{2p_2+1}, \cr
&    \qquad \text{ and } \, t > T_q,  \, \text{ for all } \, q=2p_1+1,...,2p_2 \big\},\cr
& 	{\Omega}_2(v) = \big\{(x,t) \, / \, \psi_{2p_1-1}(t) < x < \psi_{2p_2}(t),  
\, \, p_1 \le p_2, \, t<T_{2p_1-1}, \, t<T_{2p_2}, \cr
&   \qquad \text{ and } \, t > T_q,  \, \text{ for all } \, q=2p_1,...,2p_2-1 \big\},\cr
&	{\Omega}_3(v) = \text{ Closure }\big\{ (\psi_q(t),t), 
\, \text{ for all relevant values of } \, t \, \text{ and } \, q \big\}. \cr
\endaligned
$$
The set ${\Omega}_3(v)$, by construction, 
contains all of the curves $\psi_{q}$ including their end points. 
The sets ${\Omega}_1(v)$ and ${\Omega}_2(v)$ are open 
and contain regions limited by curves in ${\Omega}_3(v)$. 
These definitions take into account the fact 
that the path need not be defined for all times. 
Observe also that an arbitrary point in ${\Omega}_3(v)$
need not be a point of extremum value for $v$. 
The decomposition under 
consideration is {\it not\/} quite the obvious partition of the $(x,t)$-plane 
into regions of monotonicity for $v$. Strictly speaking, 
the sets ${\Omega}_j(v)$ may not be 
determined from the function $v$ alone.

Using Lemmas 3.2, 3.4, and 3.5, we immediately check that:

\proclaim{Lemma 3.6}  The limiting functions satisfy the properties: 
$$
\aligned
& v_{/\Omega_1(v)}(t) \, 
\text{ is non-decreasing in each subcomponent of } \, \Omega_1(v), \cr
& v_{/{\Omega}_2(v)}(t) \, 
\text{ is non-increasing in each subcomponent of } \, \Omega_2(v), 
\endaligned
\tag 3.12
$$
and
$$
w_q \, \text{ is non-decreasing if } \, q \, \text{ is even 
and non-increasing if } \, q \, \text{ is odd.}
\tag 3.13
$$
$\hfill \square$
\endproclaim

Since the scheme is consistent, conservative, and converges in the $L^1$ strong norm,
we can pass to the limit in (2.1). 
It follows that $v$ is a weak solution to (1.1). 
We note that, in the set ${\Omega}_1(v)$, 
the functions $u^h$ and, so $v$, are non-increasing. 
The Lax shock inequality holds for both $u^h$ and $v$.  
On the other hand, the cell entropy inequality (2.9) holds for $u^h$ 
in  the non-decreasing regions, i.e., in ${\Omega}_2(v)$. 
The passage to the limit in (2.9) is a classical matter. 

\proclaim{Lemma 3.7 } 
The function $v$ is a weak solution to  equation $(1.1)$ and 
satisfies 
$$
v(x-,t) \ge v(x+,t) \quad \text{ in the set } \, {\Omega}_1(v)
\tag 3.14
$$
and 
$$
\del_t U (v) + \del_x F(v) \le 0 \quad \text{ in the set } \, {\Omega}_2(v).
\tag 3.15
$$ 
$\hfill \square$
\endproclaim

The rest of this section is devoted to proving that the Lax shock inequality
holds along the paths $\psi_q$ which we will attain in Lemma 3.10. 
In a first stage, we prove: 

\proclaim{Lemma 3.8}
Along each path of extremum $\psi_q$ and for almost every $t \in (0,T_q)$,
one of the followings hold: 
$$
w_q(t) = v(\psi_q(t)-,t) \quad \text{ or } \quad  w_q(t) = v(\psi_q(t)+,t). 
\tag 3.16
$$
$\hfill \square$
\endproclaim

Roughly speaking (3.16) means that that no cusp-like layer can form in the scheme
nearby local extrema. 
The idea of the proof of Lemma 3.8 is as follows:
we are going 
to integrate the discrete form of the conservation law (2.1) on a (small) domain 
limited on one side by an approximate path of extremum, 
then we shall integrate by parts and pass to the limit as $h \to 0$. 
Finally, we shall let the domain shrink and reduce to the path itself.  
To determine the limits of the relevant boundary terms as $h \to 0$, 
we have to justify the passage to the limit in particular 
in the numerical fluxes evaluated along the approximate path. 
Lemma 3.9 below provides us with  
an a priori estimate for the oscillation of $u^h$ along the path,
which follows from the quadratic decay property (2.12).

\proclaim{Lemma 3.9} Along a path of extremum values $\psi_q^h$, we have 
$$
\beta \sum_{n=n_-}^{n=n_+} \min \bigl( |u^n_{J_q^h(n)\pm 1} - u^n_{J_q^h(n)}|,
 |u^n_{J_q^h(n)\pm 1} - u^n_{J_q^h(n)}|^2 \bigr) 
\le u^{n_-}_{J_q^h(n_-)} - u^{n_+}_{J_q^h(n_+)+1}, 
  \tag 3.17
$$
for all  $0 \le n_- \le n_+ \le N_q^h$, where
$$
\beta = \min(\alpha, 1/2). 
\tag 3.18 
$$ 
$\hfill \square$ 
\endproclaim

\proclaim{Proof of Lemma 3.8} \rm  
We will prove that, for almost every $t$ in $(0,T_q)$, 
the following three Rankine-Hugoniot like relations hold: 
$$
- \frac{d\psi_q}{dt}(t)\bigl(v(\psi_q(t)+,t) - v(\psi_q(t)-,t)\bigr)
+ f(v(\psi_q(t)+),t) - f(v(\psi_q(t)-,t)) = 0,
\tag 3.19
$$
$$
- \frac{d\psi_q}{dt}(t)\bigl(v(\psi_q(t)\pm,t) - w_q(t)\bigr)  
+ f(v(\psi_q(t)\pm,t)) - f(w_q(t)) = 0.
\tag 3.20 
$$
Since there is {\it only one\/} non-trivial pair of values that achieves 
a Rankine-Hugoniot relation for a scalar conservation law 
with a {\it strictly\/} convex flux and a given shock speed $d\psi_q(t)/dt$, 
the desired conclusion (3.16) follows immediately from (3.19)-(3.20).  

Observe that (3.19) is nothing but the standard Rankine-Hugoniot relation
since $v$ is a weak solution to (1.1)
and $\psi_q$ is Lipschitz continuous. 
For definiteness we prove (3.20) in the case of the ``$+$'' sign. 
The proof of (3.20) with  ``$-$'' sign is entirely similar.
(Actually only one of the two relations in (3.20) suffice for the present 
proof.)
 
Let $\theta(x,t)$ be a test-function having 
its support included in a neighborhood of the curve $\psi_q$ and
included in the strip $\RR\times(0,T_q)$. So for $h$ small enough 
the support of $\theta$ is included into $\RR\times(0,T_q^h)$ and all the quantities
to be considered below make sense. 
Let us set $\theta_j^n = \theta(x_j,t_n)$. 
To make use of estimate (3.17), 
it is necessary to define a ``shifted'' path $\tilde\psi_q^h$, 
to be used instead of $\psi_q^h$. 
So we consider the set of indices 
$$
\PP_q^h = \big\{ (j,n)\,  / \, j \ge J_q^h(n)+\epsilon_q^h(n) \big\},
$$
where $\epsilon_q^h(n)=0$ (respectively $\epsilon_q^h(n)=1$) if $J_q^h(n)-1$ 
(resp. $J_q^h(n)+1$) achieves the minimum in the left hand side of (3.17). 
A shifted path is defined by
$$
\tilde\psi_q^h(t) = x_{J_q^h(n)+\epsilon_q^h(n)} + \frac{t-t_n}{\tau} 
\bigl( x_{J_q^h(n+1)+\epsilon_q^h(n+1)}- x_{J_q^h(n)+\epsilon_q^h(n)} \bigr), 
\quad t \in [t_n,t_{n+1}]. 
$$
Introducing the shifts $\epsilon_q^h(n)$  does not 
modify the convergence properties of the path. It is not hard to see,
using solely the fact that the path is uniformly bounded
in Lipschitz norm, that as $ h \to 0$ 
$$
\tilde\psi_q^h  \to \psi_q  \quad W^{1,\infty} 
\text{ weak-}\star.
\tag 3.21
$$
We also set 
$$
\tilde w_q^h(t) = u_{J_q^h(n)+\epsilon_q^h(n)}^n \qquad 
\text{ for } \qquad t_n \le t < t_{n+1}. 
$$
Using (3.17) and (3.7) of $w_q^h$, it is checked that
$$
\tilde w_q^h  \to w_q  \quad L^1 \text{ strongly.} 
\tag 3.22
$$

Consider 
$$
I^h(\theta) \equiv - \sum_{(n,j) \in \PP_q^h} \bigl(u^{n+1}_j-u^n_j+
\lambda(g^n_{j+{1/2}}-g^n_{j-{1/2}})\bigr)\theta_j^n h=0, 
\tag 3.23 
$$
which vanishes identically in view of (2.1). Using summation by parts gives 
$$
\aligned
& I^h(\theta) 
= - \sum_{(n,j) \in \PP_q^h} \big(u^{n+1}_j \theta^{n+1}_j - u^n_j \theta_j^n\big) h
+\sum_n g^n_{J_q^h(n)+\epsilon_q^h(n)-{1/2}}\theta^n_{J_q^h(n)} \tau\\
& \qquad 
+ \sum_{(n,j) \in \PP_q^h} u^{n+1}_j(\theta^{n+1}_j-\theta_j^n) h + 
g^n_{j+{1/2}}(\theta^n_{j+1}-\theta_j^n) \tau\\
& \qquad  = I_1^h(\theta) + I_2^h(\theta) + I_3^h(\theta).
\endaligned
\tag 3.24
$$
The passage to the limit in $I_3^h(\theta)$ is an easy matter, 
since it has the classical form met, for instance, 
in the Lax-Wendroff theorem. We find
$$
I_3^h(\theta) \to I_3(\theta) = \iint_{\big\{x\ge\psi_q(t)\big\}}
\big(v \del_t\theta + f(v) \del_x\theta \big) dxdt. 
\tag 3.25
$$ 
To deal with $I_2^h(\theta)$, we recall that the flux
$g^n_{J_q^h(n)+\epsilon_q^h(n)-{1/2}}$ depend on two arguments so satisfies
$$
\aligned 
g^n_{J_q^h(n)+\epsilon_q^h(n)-{1/2}}  
& = f(u^n_{J_q^h(n)+\epsilon_q^h(n)}) 
+ O\big(|u^n_{J_q^h(n)-1+\epsilon_q^h(n)} - u^n_{J_q^h(n)+\epsilon_q^h(n)}|\big) \\
& = f(\tilde w_q^h(t_n))
+ O\big(|u^n_{J_q^h(n)-1+\epsilon_q^h(n)} - u^n_{J_q^h(n)+\epsilon_q^h(n)}|\big).
\endaligned 
\tag 3.26 
$$
Indeed, by construction,
the point $J_q^h(n)+\epsilon_q^h(n)-{1/2}$
is an end point of a cell achieving an extremum value:
One has $\epsilon_q^h(n) \in \{0,1\}$, so 
$J_q^h(n)+\epsilon_q^h(n)-{1/2} \in \{J_q^h(n)-1/2, J_q^h(n)+1/2\}$. 
Using (3.26), estimate (3.17), the Cauchy-Schwartz inequality, 
and finally Lebesgue convergence theorem, 
it is not hard to prove that 
$$
I_2^h(\theta) \to I_2(\theta) = \int_{\RR_+} f(w_q(t)) \theta(\psi_q(t),t) dt. 
\tag 3.27
$$

It remains to prove that
$$
I_1^h(\theta) \to I_1(\theta) = - \int_{\RR_+} 
\frac {d\psi_q}{dt}(t) w_q(t) \theta(\psi_q(t),t) \, dt. 
\tag 3.28
$$
We return to the definition of the modified path
and define $e_q^h(n+1)$ by 
$$
J_q^h(n+1)+\epsilon_q^h(n+1) = J_q^h(n)+\epsilon_q^h(n) + e_q^h(n+1).
$$
Using only (3.21), one can prove
$$
\sum_n e_q^h(n) u^n_{J_q^h(n)+\epsilon_q^h(n)}
\theta^n_{J_q^h(n)+\epsilon_q^h(n)} h \, \to \, 
- \int_{\RR_+} \frac {d\psi_q}{dt}(t)  \theta(\psi_q(t),t)  \, dt. 
\tag 3.29 
$$
We claim that
$$
e_q^h(n+1) \in \big\{-1,0,1\big\}. 
\tag 3.30
$$
Namely, using (2.11) and the definition of $\epsilon_q^h(n+1)$, 
we have either
$$
\epsilon_q^h(n)=1 \quad  
\quad J_q^h(n+1)=  J_q^h(n) \, \text{ or }  \, J_q^h(n+1)+1, 
\text{ and } e_q^h(n+1) = 1  \, \text{ or }  \, 1
$$ 
or 
$$
\epsilon_q^h(n)=1 \quad  
\quad J_q^h(n+1)=  J_q^h(n) \, \text{ or }  \, J_q^h(n+1)+1, 
\text{ and } e_q^h(n+1) = - 1  \, \text{ or }  \, 0. 
$$ 

The term $I_1^h(\theta)$ then can be rewritten in the form
$$
\align
I_1^h(\theta) 
& = \sum_{j \ge J_q^h(n)+\epsilon_q^h(n)} u^{n+1}_j \theta^{n+1}_j h - 
      \sum_{j \ge J_q^h(n)+\epsilon_q^h(n)} u^n_j \theta_j^n h\cr
& = \sum_{j \ge J_q^h(n-1)+\epsilon_q^h(n-1)} u^n_j \theta^n_j h -
      \sum_{j \ge J_q^h(n)+\epsilon_q^h(n)} u^n_j \theta_j^n h, \cr
\endalign
$$
so that 
$$
I_1^h(\theta) 
= - \sum_{e_q^h(n)=-1} u^n_{J_q^h(n)+\epsilon_q^h(n)}
 \theta^n_{J_q^h(n)+\epsilon_q^h(n)} h 
+ \sum_{e_q^h(n)=1} u^n_{J_q^h(n)+
\epsilon_q^h(n)-1} \theta^n_{J_q^h(n)+\epsilon_q^h(n)-1} h.
$$
Observe that, in the second sum above, 
$J_q^h(n)+\epsilon_q^h(n)-1 = J_q^h(n-1)+\epsilon_q^h(n-1)$, and
consider the decomposition
$I_1^h(\theta) = I_{1,1}^h(\theta) + I_{1,2}^h(\theta)$
with
$$
I_{1,1}^h(\theta) = \sum_{e_q^h(n)=1} \big(u^n_{J_q^h(n)+\epsilon_q^h(n)-1}
 \theta^n_{J_q^h(n)+\epsilon_q^h(n)-1} - u^n_{J_q^h(n)+\epsilon_q^h(n)}
 \theta^n_{J_q^h(n)+\epsilon_q^h(n)}\big) h
$$
and
$$
I_{1,2}^h(\theta) = \sum_{e_q^h(n)=1} u^n_{J_q^h(n)+\epsilon_q^h(n)}
 \theta^n_{J_q^h(n)+\epsilon_q^h(n)} h - \sum_{e_q^h(n)=-1} u^n_{J_q^h(n)+
\epsilon_q^h(n)} \theta^n_{J_q^h(n)+\epsilon_q^h(n)} h.
$$

On one hand, we have 
$$
\align
I_{1,1}^h(\theta) & \le \sum_{e_q^h(n)=1} \big|u^n_{J_q^h(n)+\epsilon_q^h(n)-1} 
\theta^n_{J_q^h(n)+\epsilon_q^h(n)-1} - u^n_{J_q^h(n)+\epsilon_q^h(n)} 
\theta^n_{J_q^h(n)+\epsilon_q^h(n)}\big| h \cr
& \le O(1) \sum_{e_q^h(n)=1} \big|u^n_{J_q^h(n)+\epsilon_q^h(n)-1} -  
u^n_{J_q^h(n)+\epsilon_q^h(n)} \big| h  \cr
& + O(1)\sum_{e_q^h(n)=1} \big|\theta^n_{J_q^h(n)+\epsilon_q^h(n)-1} - 
\theta^n_{J_q^h(n)+\epsilon_q^h(n)}\big| h,  \cr
\endalign
$$
and, in view of (3.17) and the smoothness property of $\theta$, 
$$
\align    I_{1,1}^h(\theta) 
&\le O(1) h^{1/2}\big(\sum_n |u^n_{J_q^h(n)-1} - u^n_{J_q^h(n)}|^2\big)^{1/2} 
+ O(1) h \cr
&\le O(1) \big( h+h^{1/2}\big) \le O(1) h^{1/2},
\endalign
$$
which implies 
$$
I_{1,1}^h(\theta)  \to 0  \quad \text {as } \,  h  \to  0.  
\tag 3.31
$$
The expression for $I_{1,2}^h(\theta)$ can be simplified, namely
$$
I_{1,2}^h(\theta) =  - \sum_n e_q^h(n) u^n_{J_q^h(n)+\epsilon_q^h(n)}
 \theta^n_{J_q^h(n)+\epsilon_q^h(n)} h.
$$
Using (3.29), we find that
$$
I_{1,2}^h(\theta)  \to  I_1(\theta)   \quad \text {as } \,  h  \to  0.
\tag 3.32
$$

In view of (3.25), (3.27), and (3.28),  we conclude that 
$$
I_1(\theta) + I_2(\theta) + I_3(\theta) = 0. 
\tag 3.33
$$
Finally, using in (3.33) a sequence of test-functions $\theta$, 
whose supports shrink and concentrate on the curve $\psi_q$, 
the desired Rankine-Hugoniot relation (3.20) 
with the ``$+$'' sign follows at the limit. 
This completes the proof of Lemma 3.8.
$\hfill \square$
\endproclaim

\proclaim{Proof of Lemma 3.9} \rm  
For definiteness, we assume that $u^n_{J_q^h(n)}$ is a maximum value and that:
$$
 \min \bigl( u^n_{J_q^h(n)} - u^n_{J_q^h(n)\pm 1}\bigr) 
 =  u^n_{J_q^h(n)} - u^n_{J_q^h(n)-1}.
$$
The other cases are treated similarly. To simplify the notation, set $j_*=J_q^h(n)$. 
By the uniform decay property (2.12), we have
$$
u^n_{j_*} - u^{n+1}_{j_*} \ge \alpha (u^n_j{_*}-u^n_{j{_*}-1})^2. 
\tag 3.34 
$$
By the strong maximum principle (2.11), we have
$$
u^{n+1}_{j_*-1} -  u^{n}_{j_*-1} \le \frac 1 2 (u^n_{j_*} - u^{n}_{j_*-1}).
\tag 3.35
$$
{}From (3.35), we deduce
$$
\align
u^{n+1}_{j_*-1} & \le  u^{n}_{j_*-1} + \frac 1 2 (u^n_{j_*} - u^{n}_{j_*-1})\cr
& =  u^{n}_{j_*} - \frac 1 2 (u^n_{j_*} - u^{n}_{j_*-1})\cr
& \le u^{n}_{j_*} - \frac 1 2 \min\big( (u^n_{j_*} - u^{n}_{j_*-1}), (u^n_{j_*} - u^{n}_{j_*-1})^2\big)\cr
\tag 3.36
\endalign
$$
Note that $u^n_{j_*} - u^{n}_{j_*-1}$ might be either $\le 1$ or $\ge 1$. Then we get
$$
u^{n+1}_{j_*-1} \le u^{n}_{j_*} - \beta \min\big((u^n_{j_*} - u^{n}_{j_*-1}), 
(u^n_{j_*} - u^{n}_{j_*-1})^2\big).
\tag 3.37
$$
On the other hand, (3.34) can be written in the form
$$
u^{n+1}_{j_*} \le u^n_{j_*} - \alpha (u^n_j{_*}-u^n_{j{_*}-1})^2. 
\tag 3.38
$$
Moreover, using (2.11) again, we obtain:
$$
\aligned
u^{n+1}_{j_*+1} & \le u^n_{j_*+1} + \frac 1 2  \big(u^n_{j_*} - u^{n}_{j_*+1}\big)\\
& \le u^n_{j_*} + \frac 1 2  \big(u^n_{j_*} - u^{n}_{j_*-1}\big).
\endaligned
$$
Hence, in view of (3.36):
$$
u^{n+1}_{j_*+1} \le u^{n}_{j_*} - \beta \min\big((u^n_{j_*} - u^{n}_{j_*-1}), 
(u^n_{j_*} - u^{n}_{j_*-1})^2\big).
\tag 3.39
$$
It follows from (3.37)--(3.39) that
$$
\max\big(u^{n+1}_{j_*-1},u^{n+1}_{j_*},u^{n+1}_{j_*+1}\big) \le u^{n}_{j_*} - 
\beta \min\big( (u^n_{j_*} - u^{n}_{j_*-1}), (u^n_{j_*} - u^{n}_{j_*-1})^2\big),
\tag 3.40
$$
and thus
$$
u^{n+1}_{J_q^h(n+1)} \le u^{n}_{j_*} - \beta \min\big( (u^n_{j_*} - u^{n}_{j_*-1}), 
(u^n_{j_*} - u^{n}_{j_*-1})^2\big),
\tag 3.41
$$
since $u^{n+1}_{J_q^h(n+1)}$ 
by definition achieves the maximum in (3.20). Finally, we have proved
$$
u^{n}_{J_q^h(n)} - u^{n+1}_{J_q^h(n+1)} \ge \beta \min\big((u^n_{j_*} -
 u^{n}_{j_*-1}), (u^n_{j_*} - u^{n}_{j_*-1})^2\big).
\tag 3.42
$$
By assumption, $u^n_{j_*} - u^{n}_{j_*-1} \le u^n_{j_*} - 
u^{n}_{j_*}$. Thus (3.42) gives (3.17) after summation w.r.t. $n$.
$\hfill \square$
\endproclaim

\proclaim{Lemma 3.10} We have 
$$
v(x-,t) \ge v(x+,t) \quad \text{ in the set } \, {\Omega}_3(v).
\tag 3.43 
$$
$\hfill \square$
\endproclaim

The proof is based on the fact that one of the entropy criteria is 
satisfied on each side of a path: the Lax shock inequality
in the non-increasing side, and the cell entropy
inequality in the non-decreasing one. 

\proclaim{Proof}  \rm 
We claim that, along any path of extremum $\psi_q$,
$$
v(\psi_q(t)-,t) \ge v(\psi_q(t)+,t) \qquad \text{ for a.e. } t \in (0,T_q).
\tag 3.44
$$
We use the notation introduced in the proof of Lemma 3.8. 
A new difficulty arises: several paths may accumulate in a region.
Lemma 3.8 was concerned with the discrete conservation laws (2.1) 
which holds in both the non-increasing and non-decreasing regions.
For the entropy consistency,
we do not use the same criterion, and this complicates the proof.

To begin with, consider a path $\psi_q$ and a point 
$(\psi_q(t_0),t_0)$, that is supposed to be 
an ``isolated'' point of change of monotonicity, 
in the sense that: $\psi_{q-1}(t_0)<\psi_q(t_0)<\psi_{q+1}(t_0)$. 
By continuity, these inequalities then hold with $t_0$ replaced by any $t$ 
lying in a small neighborhood of $t_0$. 
For definiteness, we also suppose that $\psi_q$ is a path of minimum values. 
We later analyze the case that 
two or more paths of extrema accumulate 
in a neighborhood of $(\psi_q(t_0),t_0)$. 

Let $\theta$ be a non-negative test-function of the two variables $(x,t)$ 
having its support included in a small neighborhood of $(\psi_q(t_0),t_0)$.
We can always assume that   
$u^h$ is non-increasing on the left side of the curve $\psi_q^h$, 
and non-decreasing on the right side. 
Using the notation introduced in the proof of Lemma 3.8, 
we aim at passing to the limit in 
$$
I^h(\theta) = -\sum_{(j,n) \in \PP_q^h} \big(U^{n+1}_j-U^n_j+
\lambda \, (G^n_{j+{1/2}}-G^n_{j-{1/2}})\big) \, \theta_j^n h \ge 0.
$$
Note that $I^h(\theta)$ is non-positive according to 
the cell entropy inequality (2.9)
and since the $u^h$'s are non-decreasing on the right side.

Integrating by parts in $I^h(\theta)$ gives 
$$
\aligned
& I^h(\theta) = - \sum_{(j,n) \in \PP_q^h} 
\big(U^{n+1}_j \theta^{n+1}_j - U^n_j \theta_j^n\big) h
+  \sum_n G^n_{J_q^h(n)+\epsilon_q^h(n)-{1/2}} \, 
\theta^n_{J_q^h(n)+\epsilon_q^h(n)-{1/2}}  \, \tau\\
& \qquad 
+ \sum_{(j,n) \in \PP_q^h} U^{n+1}_j(\theta^{n+1}_j-\theta_j^n) h +
G^n_{j+{1/2}} \, (\theta^n_{j+1}-\theta_j^n) \tau\\
& \qquad  = I_1^h(\theta) + I_2^h(\theta) + I_3^h(\theta).
\endaligned
\tag 3.45
$$
The passage to the limit in the term $I_3^h(\theta)$ is a classical matter. 
The treatment of $I_2^h(\theta)$ and  $I_1^h(\theta)$ 
is similar to what was done to prove (3.27) and (3.28), respectively. 
Therefore we have 
$$
I_1^h(\theta) \to I_1(\theta) = - \int_{\RR_+} \frac {d\psi_q}{dt}(t) \, 
U(w_q(t)) \, \theta(\psi_q(t),t) \, dt. 
\tag 3.46
$$
$$
I_2^h(\theta) \to I_2(\theta) = \int_{\RR_+} F(w_q(t)) \, 
\theta(\psi_q(t),t) \, dt,
\tag 3.47
$$
$$
I_3^h(\theta) \to I_3(\theta) = \iint_{\big\{x\ge\psi_q(t)\big\}} 
\big(U(v) \del_t\theta + F(v) \del_x\theta \big) dxdt. 
\tag 3.48
$$
It follows from (3.46)--(3.48) that
$$
I_1(\theta) + I_2(\theta) + I_3(\theta) \ge 0. 
\tag 3.49
$$
Finally, using a sequence of test-functions whose supports
shrink and concentrate on the curve $\psi_q$, we deduce from (3.49)
that the entropy dissipation is non-positive along the path, i.e.
$$
- \frac{d\psi_q}{dt}(t)\bigl(U\big(v(\psi_q(t)+,t)\big) - 
U\big(w_q(t)\big)\bigr) + 
F(v(\psi_q(t)+,t)) - F(w_q(t)) \le 0. 
$$
Combined with (3.20), this inequality is equivalent 
to  
$$
v(\psi_q(t)+,t)) \le  w_q(t), 
$$
which yields the desired inequality (3.44). 

Consider next the case that several paths accumulate 
in the neighborhood of $(\psi_q(t_0),t_0)$. 
In view of (1.4), a finite number of paths only
can accumulate at a given point. 
For definiteness we suppose that \par
-- the point $(\psi_q(t_0),t_0)$ is a point of minimum values for $u^h$; \par
-- the curves $\psi_{q+1}$ and $\psi_{q+2}$ coincide with $\psi_q$ 
in a neighborhood of $t_0$; \par 
-- and we have $\psi_{q-1}<\psi_q$ and $\psi_{q+2}<\psi_{q+3}$ 
in a neighborhood of $t_0$.  \par 
\noindent Suppose that, for instance, $w_q(t) \le w_{q+2}(t)$. 
The other cases are treated similarly. 
In this situation, we are going to prove that
$$
v(\psi_q(t)-,t) \ge  w_q(t)= w_{q+1}(t) \ge w_{q+2}(t)= v(\psi_q(t)+,t), 
\tag 3.50
$$
at those points $t$ near $t_0$ where
$\psi_q(t)=\psi_{q+1}(t)=\psi_{q+2}(t)$. Of course, (3.41) is a much 
stronger statement than (3.34). 
By definition of the paths of extrema, we have 
$$
w_q(t) \le  w_{q+1}(t) \quad \text {and} \quad w_{q+1}(t) \ge w_{q+2}(t),  
$$
$$
v(\psi_q(t)-,t) \ge  w_q(t) \quad \text {and} \quad w_{q+2}(t) \le v(\psi_q(t)+,t).
$$
Lemma 3.8 shows that
$$
w_q(t), w_{q+1}(t), \, \text {and} \, w_{q+1}(t) \in 
\big\{ v(\psi_q(t)-,t), v(\psi_q(t)+,t)\big\}. 
$$
Thus, in order to get (3.50), 
it is sufficient to check the following two inequalities 
$$
w_{q+2}(t) \ge v(\psi_q(t)+,t), 
\tag 3.51 
$$
$$
w_q(t) \ge w_{q+1}(t). 
\tag 3.52 
$$

On one hand, 
the argument used in the first part of the present proof applies directly to the 
path $\psi_{q+2}(t)$ and the region located to the right of this curve,
since $\psi_{q+2}(t)<\psi_{q+3}(t)$ in a neighborhood of $t_0$. 
As a consequence, we obtain
$w_{q+2}(t) \ge v(\psi_{q+2}(t)+,t)$, 
which is exactly (3.51), since $\psi_{q+2}= \psi_q$. 

On the other hand, to prove (3.52), let $\PP_{q,q+1}^h$ be 
the (small) region limited by the curves $\psi_q^h$ and $\psi_{q+1}^h$ for
$t$ belonging to a small neighborhood of $t_0$. Specifically, $\PP_{q,q+1}^h$ 
is a set of indices of the form $(j,n)$ defined along the lines of 
the proof of Lemma 3.8. In particular, both paths are modified according 
to estimate (3.17), as explained before. Consider 
$$
I^h(\theta) = - \sum_{(j,n) \in \PP_{q,q+1}^h} \big(U^{n+1}_j-U^n_j+
\lambda \, (G^n_{j+{1/2}}-G^n_{j-{1/2}})\big) \, \theta_j^n h \ge 0.
\tag 3.53 
$$
Note that $I^h(\theta)$ is non-negative and
$$
\aligned
& I^h(\theta)  = -  \sum_{(j,n) \in \PP_{q,q+1}^h} 
\big(U^{n+1}_j \theta^{n+1}_j - U^n_j \theta_j^n\big) h \\ 
& + \sum_n G^n_{J_q^h(n)+\epsilon_q^h(n)-{1/2}} \, 
\theta^n_{J_q^h(n)+\epsilon_q^h(n)-{1/2}}  \, \tau
- \sum_n G^n_{J_{q+1}^h(n)+\epsilon_{q+1}^h(n)-{1/2}} \, 
\theta^n_{J_{q+1}^h(n)+\epsilon_{q+1}^h(n)-{1/2}}  \, \tau \\ 
& + \sum_{(j,n) \in \PP_{q,q+1}^h} U^{n+1}_j(\theta^{n+1}_j-\theta_j^n) h +
G^n_{j+{1/2}} \, (\theta^n_{j+1}-\theta_j^n) \tau\\
& = I_1^h(\theta) + I_{2,q}^h(\theta) - I_{2,q+1}^h(\theta) + I_3^h(\theta).
\endaligned 
$$
Using the technique developed for the proof of Lemma 3.8, we get
$$
\aligned 
I_1^h(\theta) \to  I_1(\theta) 
=  - & \int_{\RR_+} \frac {d\psi_q}{dt}(t) \, U(w_q(t)) \, \theta(\psi_q(t),t) \, dt \\
   +  &  \int_{\RR_+} \frac {d\psi_{q+1}}{dt}(t) \, U(w_{q+1}(t)) \, 
\theta(\psi_{q+1}(t),t) \, dt,
\endaligned
$$
$$
\aligned 
&  I_{2,q}^h(\theta) \to  I_{2,q}(\theta) = 
\int_{\RR_+} F(w_q(t)) \, \theta(\psi_q(t),t) dt,\\
&  I_{2,q+1}^h(\theta) \to  I_{2,q+1}(\theta) = 
  - \int_{\RR_+} F(w_{q+1}(t)) \, \theta(\psi_{q+1}(t),t) \, dt,\\
&  I_3^h(\theta) \to 0. 
\endaligned 
$$
It follows that
$$
I_1(\theta) + I_{2,q}(\theta) + I_{2,q+1}(\theta) \ge 0, 
\tag 3.54 
$$
which, since  $\psi_q= \psi_{q+1}$ near $t_0$,
is equivalent to the jump condition 
$$
-\frac {d\psi_q}{dt}(t) \,\bigl( U(w_{q+1}(t)) - U(w_{q}(t))\bigr) 
+ F(w_{q+1}(t)) - F(w_q(t)) \le 0, 
$$
which gives (3.52). 
This completes the proof of Lemma 3.10. 
$\hfill \square$
\endproclaim

It is a classical matter to check that the initial condition (1.2) is satisfied
by using (2.2) and the uniform BV bound. 
Since the function $v(t)$ has bounded total variation, it admits left and right 
traces at each point and (1.5) and (1.6) are known to be equivalent 
at a point of discontinuity. Therefore the
following result follows from Volpert's proof 
in \cite{\refVolpert}.

\proclaim{Lemma 3.11} Supose $v$ is a 
function of bounded variation and a weak solution to 
the conservation law $(1.1)$ and satisfies the initial condition $(1.2)$, 
and the inequalities $(3.14)$, $(3.15)$, and 
$(3.34)$ where $\Omega_1(v) \cup \Omega_2(v)\cup\Omega_3(v) = \RR\times\RR_+$. 
Then $v$ is the unique entropy solution to $(1.1)$-$(1.2)$. 
$\hfill \square$
\endproclaim

The proof of Theorem 2.3 is now complete. 

The limiting paths $\psi_q$ associated with the scheme determine a decomposition 
of the plane into non-increasing/non-decreasing regions for the exact solution $u$. 
Such a decomposition is not unique, in general. Consider the decomposition 
found in the Appendix for the function $v$ and the corresponding paths $\varphi_p$. 
When $v$ is not constant in any neighborhood of an extremum path $\varphi_p$, 
the path is unique and must coincide with one of the path $\psi_q$.  
When $v$ is constant in the neighborhood of a path $\varphi_p$, then 
the path may be arbitrarily modified and it may happen that no limiting path $\psi_q$ 
coincide with $\varphi_p$. 
 
%---------------------------------------------------------------------
\heading{4. Application to the MUSCL Scheme}
\endheading 

The purpose of this section
is to apply Theorem 2.3 to van Leer's MUSCL scheme
(for Monotone Upstream Scheme for Conservation Laws); 
cf.~\cite{\refvanLeerone, \refvanLeertwo}. 
This section also provides a proof of estimate (2.12) stated in Proposition 2.1,
a new property of the Godunov scheme which does also hold for the MUSCL scheme.

It is convenient to formulate (2.1) 
in terms of the incremental coefficients $C^{\pm, n}_{j+1/2}$ defined by
$$
C^{+,n}_{j+1/2} = - \, \lambda \, \frac {g_{j+1/2}^n - f(u^n_j)}{u^n_{j+1}-u^n_j}, \quad 
C^{-,n}_{j-1/2} =  \, \lambda \, \frac {f(u^n_j)-g^n_{j-1/2}}{u^n_j-u^n_{j-1}},
\tag 4.1 
$$
so that 
$$
u^{n+1}_j = u^n_j + C^{+,n}_{j+1/2} \big(u^n_{j+1}-u^n_j\big)   
                          + C^{-,n}_{j-1/2}\big(u^n_{j-1}-u^n_j).
\tag 4.2 
$$
The numerical viscosity coefficient (Cf.~Tadmor \cite{\refTadmorone}) 
being defined by  
$$
Q_{j+1/2}^n = C^{+,n}_{j+1/2} + C^{-,n}_{j+1/2},  
\tag 4.3 
$$   
the viscous form of the scheme is 
$$
u^{n+1}_j = u^n_j - {\lambda 
\over 2} \, \bigl(f(u^n_{j+1}) - f(u^n_{j-1})\bigr) 
+ {1 \over 2} Q^{n}_{j+1/2} \big(u^n_{j+1}-u^n_j\big)   
- {1 \over 2} Q^{n}_{j-1/2}\big(u^n_{j-1}-u^n_j).
$$

\proclaim{Proposition 4.1 } 
The scheme $(2.1)$ satisfies the local maximum principle $(2.11)$
provided 
$$
C^{+,n}_{j+1/2} \ge 0, \quad 
C^{-,n}_{j-1/2} \ge 0, \quad \text{ and } \quad  
C^{+,n}_{j+1/2}+C^{-,n}_{j-1/2} \le \frac 1 2. 
\tag 4.4 
$$
$\hfill \square$. 
\endproclaim

A sufficient condition for (4.4) to hold is 
$$
C^{\pm,n}_{j+1/2} \ge 0 \quad \text{ and } \quad Q^n_{j+1/2} \le 1/4. 
\tag 4.5  
$$
Namely, if (4.5) holds, then 
$0 \le C^{+,n}_{j+1/2}\le 1/4$ and $0 \le C^{-,n}_{j+1/2} \le 1/4$,
so that (4.4) is satisfied. 
In particular, the Godunov and Engquist-Osher schemes 
satisfy (4.4) under the CFL condition (2.4).  When the numerical flux
in independent of $\lambda$, the second inequality 
in (4.5) is always satisfied provided $\lambda$ is small enough.

The Lax-Friedrichs type scheme have a constant numerical viscosity 
$Q^n_{j+1/2} \equiv Q$. For the original Lax-Friedrichs scheme 
$Q=1$. Proposition 4.1 applies provided $ Q \leq 1/4$. Observe that 
the monotonicity property {\it does fail\/} when $Q \in(2/3, 1]$: 
take for instance $f \equiv 0$ and 
$u_j^n \equiv 0$ for all $j \ne 0$ but $u_0^n >0$.  
This initial data has one maximum point, and at the next time step 
$$
u_j^{n+1} = {Q\over 2} u^n_{j+1} + (1-Q) u^n_j + {Q\over 2} u^n_{j-1}. 
$$
admits two maximum points $j=-1$ and $j=1$. 
A related observation was made by Tadmor in \cite{\refTadmorone}: 
for $Q \leq 1/2$, better properties can be obtained for the Lax-Friedrichs 
scheme.

\proclaim{Proof of Proposition 4.1} \rm 
Inequalities (2.11) can be written in terms of the 
incremental coefficients, namely 
$$
\aligned 
\frac 1 2 \min\big(\delta^n_{j+1}, \delta^n_j, \delta^n_{j-1}\big) 
& \le  C^{+,n}_{j+1/2} \delta^n_{j+1}  
+ \big(\frac 1 2 -C^{+,n}_{j+1/2}-C^{-,n}_{j-1/2}\big) \delta^n_j
+ C^{-,n}_{j-1/2}\delta^n_{j-1} \\ 
& \le \frac 1 2 \max\big(\delta^n_{j+1}, \delta^n_j, \delta^n_{j-1}\big)
\endaligned 
\tag 4.6 
$$
with
$$
\delta^n_{j+1} = u^n_{j+1} - u^n_j, \quad 
\delta^n_j = 0, \quad
\delta^n_{j-1} = u^n_{j-1} - u^n_j. 
$$  
If (4.4) holds, then   
$$
2 \, C^{+,n}_{j+1/2} \delta^n_{j+1}  
+ \big(1 - 2 \, C^{+,n}_{j+1/2 }-  2 \, C^{-,n}_{j-1/2}\big) \delta^n_j
+ 2 \, C^{-,n}_{j-1/2}\delta^n_{j-1}) 
$$
is a convex combination of the $\delta_j$'s. So
(4.6) and therefore (2.11) follows. 
$\hfill \square$
\endproclaim

%----------------------------------------------

We now introduce the van Leer's scheme, composed 
of a reconstruction step based on the min-mod limitor and a resolution step
based on the Godunov solver. 
We use the notation introduced in Section 2. For simplicity in the presentation, 
we normalize the flux to satisfy $f(0)= f'(0) = 0$. 
{}From the approximation $\big\{u_j^n\big\}$ at the time $t=t_n$,
we construct a piecewise affine function 
$$
\tilde u^n_j(x) = u_j^n + s_j^n (x-x_j)/h \quad \text { for } \, 
x \in (x_{j-1/2}, x_{j+1/2}), 
\tag 4.7 
$$
where the slope $s_j^n$ is 
$$
s_j^n = \minmod \big(u_j-u_{j-1}, (u_{j+1} - u_{j-1})/2, u_{j+1} - u_j\big) 
\tag 4.8 
$$
with 
$$ 
\minmod (a,b,c) = 
\cases & \min (a,b,c) \quad \text{ if } \, a>0, b>0, \, \text{ and } \, c>0, \\
           & \max (a,b,c) \quad \text{ if } \, a<0, b<0, \, \text{ and } \, c<0. \\
           & 0  \quad \quad     \qquad \text{ in all other cases.} 
\endcases
\tag 4.9 
$$

We introduce the notation
$$
u^n_{j+1/2-} = u_j^n + s_j^n /2 \quad \text{ and } \quad 
u^n_{j+1/2+} = u_{j+1}^n - s_{j+1}^n /2. 
$$
Then the solution is up-dated with (2.1) where the numerical flux is 
defined depending upon the values of the reconstruction at the interfaces. 

\roster 
\item 
If either $0 \leq u^n_{j} \leq u^n_{j+1} $ or $0 \leq u^n_{j+1} \geq u^n_j$, 
then the numerical flux is defined 
by using the characteristic line traced backward
from the point $(x_{j+1/2}, t_{n+1/2})$. Since the latter 
has a positive slope, we set  
$$
g^n_{j+1/2} = f(u^n_{j+1/2-}) + f'(u^n_{j+1/2-}) \, (v- u^n_{j+1/2-}) 
\tag 4.10a 
$$
with $v$ solving 
$$
u^n_{j+1/2-} = v + \frac \lambda 2 f'(v) s_j^n.
\tag 4.10b 
$$ 
\item 
If either $u^n_{j} \leq u^n_{j+1} \leq 0$, or 
$u^n_{j+1/} \leq u^n_{j} \leq 0$, 
then the backward characteristic has a negative slope and we set
$$
g^n_{j+1/2} = f(u^n_{j+1/2+}) + f'(u^n_{j+1/2+}) \, (v - u^n_{j+1/2+}) 
\tag 4.11a 
$$
with $v$ solving 
$$
u^n_{j+1/2+} = v + \frac \lambda 2 f'(v) s_{j+1}^n.
\tag 4.11b 
$$
\item In all other cases we set 
$$
g^n_{j+1/2} = f(0) = 0.
\tag 4.12 
$$
\endroster 
Equations (4.10b) and (4.11b) can be solved explicitly for the Burgers equation 
since then $f'(u) = u$ is linear. Observe that the scheme reduces to first order 
at sonic points and extrema. 

The main result of this section is:

\proclaim{Theorem 4.2 } 
For $\lambda$ small enough, the MUSCL method defined by $(4.7)$--$(4.12)$ 
is a generalized monotone scheme in the sense of Definition $2.2$.
When $(1.3)$-$(1.4)$ hold,
 the scheme converges in the strong $L^1$ topology 
to the unique entropy solution of $(1.1)$-$(1.2)$. 
$\hfill \square$
\endproclaim

It would be interesting to extend Theorem 4.2 to higher-order methods
such as the Woodward-Collela's P.P.M. scheme.

\proclaim{Proof of Theorem 4.2 } \rm
We have to check that the scheme satisfies the three conditions 
in Definition 2.2. We always assume that $\lambda$ is, at least,  
less or equal to $1/4$. 
\vskip.1cm

\noindent{ \bf Step 1: } Local maximum principle.  

Estimate (2.11) is easily obtained by applying Proposition 4.1 
and relying on the convexity of the flux function $f$. 
We omit the details. 
\vskip.1cm 

\noindent{ \bf Step 2: } Cell entropy inequality. 

Consider a region where the sequence $\big\{u^n_j \big\}$ is 
non-decreasing. We will use the entropy pair $(U,F)$ with 
$U(u) = u^2/2$ and $F'(u) = u f'(u)$. Define 
the numerical entropy flux by
$$
G^n_{j+1/2} = F(u^n_{j+1/2-}) + U'(u^n_{j+1/2-}) \, (f(v) - f(u^n_{j+1/2-})), 
\tag 4.13i
$$
$$
G^n_{j+1/2} = F(u^n_{j+1/2+}) + U'(u^n_{j+1/2+}) \, (f(v) - f(u^n_{j+1/2+})), 
\tag 4.13ii 
$$
and 
$$
G^n_{j+1/2} = 0
\tag 4.13iii
$$
in Cases (1), (2), and (3), respectively. 
Inequality (2.9) is checked by direct calculation, 
for $\lambda$ small enough. Observe that Case (3) is obvious 
since our scheme then reduces to a first order, entropy consistent scheme. 

For definiteness we treat Case (1), i.e.~$f'>0$ for the 
values of $u$ under consideration. 
We view the left hand side of (2.9) as a function of 
$w=u_{j-1/2-}$, $u=u_j^n$, $v=u_{j+1/2-}$ and the value $\tilde w$ 
defined as 
$$
w = \tilde w + {\lambda \over 2} f'(\tilde w) t, 
$$ 
where $t$ stands for the slope in the cell $j-1$.
Introduce also $\tilde v$ by 
$$ 
v = \tilde v + {\lambda \over 2} f'(\tilde v) s 
$$
with $s= 2 (v-u)$. 
Since the approximate solution is non-increasing, we have $\tilde w \leq 
w \leq 2 u - v \leq u \leq v \leq \tilde v$. Set 
$$
\aligned 
& \Omega(\tilde w, w, u,v;\lambda) \\ 
& = U(\bar u) - U(u) + \lambda \bigl[F(v) + U'(v)(f(\tilde v) - f(v)) 
- F(w) - U'(w)(f(\tilde w) - f(w))\bigr], 
\endaligned 
$$
and 
$$
\bar u = u - \lambda \bigl[f(v) + f'(v)(\tilde v - v) 
- f(w) - f'(w)(\tilde w - w)\bigr]. 
$$

Observe that 
$$
\del_{\tilde w} \Omega 
= U'(\bar u) \lambda f'(w) - \lambda U'(w) \, f'(\tilde w), 
$$
and 
$$
\del^2_{\tilde w} \Omega 
= U''(\bar u) \lambda^2 f'(w) - \lambda U'(w) \, f''(\tilde w) \leq - C \lambda \, |w| 
$$
for $w>0$ and $\lambda$ small enough. 
Therefore $\Omega$ is a concave function in $\tilde w$ and 
$$
\Omega(\tilde w, w, u,v;\lambda) 
\leq 
\Omega(w, w, u,v;\lambda) - (w-\tilde w) \del_{\tilde w}\Omega(w, w, u,v;\lambda) 
 - C \lambda |w| \, |\tilde w - w|^2. 
$$
But 
$$
\aligned 
\del_{\tilde w} \Omega (w, w, u,v;\lambda)
& = \lambda f'(w) \bigl( U'(\bar u) -  U'(w) \bigr) \\ 
& = U''(\xi) \, \lambda f'(w) 
\bigl[u-w - \lambda \bigl(f(v) + f'(v)(\tilde v - v) - f(w) \bigr)\bigr] \\ 
& \leq C \lambda \, |u-w| 
\endaligned 
$$ 
for some $\xi>0$ and $\lambda$ small enough. 
This proves that 
$$
\aligned 
\Omega(\tilde w, w, u,v;\lambda) 
& \leq 
\Omega(w, w, u,v;\lambda) - C \lambda \bigl( |\tilde w - w| |u-w|
+ |w| \, |\tilde w - w|^2\bigr)\\ 
& \leq 
\Omega(w, w, u,v;\lambda), 
\endaligned 
$$ 
and we now simply use the notation $\Omega(w, u,v;\lambda)$.

Taylor expanding $\Omega$ with respect to $\lambda$ shows 
that the dominant term is the first order 
coefficient in $\lambda$ given by 
$$
\hat \Omega_1 (w,u,v) \equiv -U'(u) \bigl(f(v) - f(w)\bigr) + F(v) - F(w), 
$$
in which $w \leq 2 u - v \leq u \leq v$. Since 
$$
\del_w \hat \Omega_1(w,u,v) = (U'(u) - U'(w)\bigr) \, f'(w)
 \geq C \, |w| \, |u-w|, 
$$
we have 
$$
\hat \Omega_1(w,u,v) \leq \hat \Omega_1(2u-v,u,v) - 
C \, |w| \int_{w}^{2u-v} (u-z) \, dz, 
$$ 
so 
$$ 
\hat \Omega_1(w,u,v) \leq \hat \Omega_1(2u-v,u,v) 
- C' \, |w| \, |2u-v -w | \, |u-w|. 
$$

It remains to study 
$\hat \Omega_1(2u-v,u,v) = \tilde \Omega_1(u,v)$ with $u \leq v$. 
We find 
$$
\aligned 
\del_v \tilde \Omega_1(u,v) & = \bigl(U'(v)-U'(u)\bigr) f'(v) 
+ \bigl(U'(2u-v) - U'(u)\bigr) f'(2u-v)\\ 
& = (v-u) \, \bigl(f'(v) - f'(2u-v)\bigr) \leq - C \, |u-v|^2. 
\endaligned 
$$ 
It follows that $\tilde \Omega(u,v)$ is a non-increasing function 
of $v$ for all $v \geq u$, and since it vanishes for $v=u$, 
$$
\tilde \Omega_1(u,v) \leq -C' \, |u-v|^3. 
$$ 

This proves that the {\it first order\/} term in $\lambda$ in the 
expansion of the function $\Omega$ is negative. 

The same arguments are now applied to the function $\Omega(\lambda)$ 
directly. We have
$$
\aligned
  \del_w \Omega_(w,u,v;\lambda) 
 & = (U'(\bar u) - U'(w)) \, \lambda f'(w)
\\
 & = \lambda f'(w)U''(\xi)(u-w-\lambda(f(v)+f'(v)(\tilde v-v) - f(w)))
\\
 & = \lambda f'(w)U''(\xi)(u-w-\lambda O(1)(u-w))
\\ 
 & \geq \lambda C \, |w| \, |u-w| \,.
\endaligned
$$
Therefore, 
$$ 
  \Omega(w,u,v;\lambda) \leq \Omega(2u-v,u,v;\lambda) 
  - \lambda C' \, |w| \, |2u-v -w | \, |u-w|. 
$$
Denote $\hat\Omega(u,v;\lambda)= \Omega(2u-v,u,v;\lambda)$
$$
 \hat \Omega(u,v;\lambda) = 
U(\bar u) - U(u) + \lambda 
\bigl[F(v) + U'(v)(f(\tilde v) - f(v))- F(2u-v)\bigr], 
$$
where 
$$
  \bar u = u - \lambda \bigl[f(v) + f'(v)(\tilde v - v) - f(2u-v)\bigr]. 
$$
We easily compute that
$$
 {\del_v \hat \Omega(u,v;\lambda)\over \lambda}
  = f'(v) \bigl(U'(v)-U'(\bar u)\bigr) + f'(2u-v) 
  \bigl(U'(2u-v) - U'(\bar u)\bigr) + \lambda A(u,v;\lambda)
$$ 
with 
$$
 |A(u,v;\lambda)| \leq C \, |u-v|^2. 
$$ 
This establishes the desired conclusion for $\lambda$ small enough.
\vskip.1cm

\noindent{ \bf Step 3: } Quadratic decay property.

Near a local extremum, the MUSCL scheme essentially reduces to 
the Godunov scheme. So it is enough to check the quadratic decay property 
(2.12) for the Godunov scheme. This can be done from the explicit formula
(2.6). 

The simplest situation is obtained with the Godunov scheme and when $f'$
has a sign, say is positive. Assume $u_j^n$ is a local maximum. We have 
$$
u^{n+1}_j = u^n_j - \lambda \bigl( f(u^n_j) - f(u^n_{j-1})\bigr), 
$$
thus
$$
\aligned 
u^n_j - u^{n+1}_j  
& = \lambda \bigl( f(u^n_j) - f(u^n_{j-1})\bigr) \\ 
& \ge 
\lambda f'(u^n_{j-1}) \, \bigl(u^n_j - u^n_{j-1}\bigr) 
+ \lambda  \bigl(\inf f''/2\bigr) \, \bigl(u^n_j - u^n_{j-1}\bigr)^2 \\
&  \ge 
\lambda  \bigl(\inf f''/2\bigr) \, \min_\pm \bigl(u^n_j - u^n_{j\pm 1}\bigr)^2.\\
\endaligned
$$
This establishes (2.11) when $f'>0$.

It remains to treat the sonic case where $f'$ has no definite sign. 
We will rely on the following technical remark. 
Given two points such that 
$$
u_- < 0 < u_+, \qquad f(u_-) = f(u_+),
$$
there exist $c_1, c_2>0$ (independent of $u_\pm$) such that 
$$
c_1 \, |u_-| \leq |u_+| \leq c_2 \, |u_-|. 
$$ 

Consider the case $u^n_{j-1} < 0 < u^n_j$, and use Osher's formula
for the Riemann problem, we have
$$
\aligned 
u^n_j - u^{n+1}_j 
& = \lambda \, (\max_{(u^n_{j+1}, u^n_j)}  f - \min_{(u^n_{j-1}, u^n_j)} f)  \\ 
& = \lambda \big(\max_{(u^n_{j+1}, u^n_j)}  f - f(0)) \\
& = \lambda  \bigl(\inf f''/2\bigr) \, (\max_{|u^n_{j+1}|, |u^n_j|} f)^2 \\ 
&\ge c |u^n_{j+1}-u^n_j|^2. 
\endaligned 
$$ 

Consider next the case $0< u^n_{j-1} < u^n_j$, then
$$
\aligned 
u^n_j - u^{n+1}_j 
 & = \lambda \big(\max_{(u^n_{j+1}, u^n_j)}  f - f(u^n_j)) \\
& \ge \lambda  \bigl(\inf f''/2\bigr) 
  \min(|u^n_{j+1}-u^n_j|^2, |u^n_{j}-u^n_{j-1}|^2)
\endaligned 
$$ 
This completes the proof of Theorem 4.2. $\hfill \square$
\endproclaim
%----------------------------------------------------------------------------

\proclaim{Acknowledgments} \rm
Most of this work was done in 1992 while P.G.L. and J.G.L. were Courant instructors
at the Courant Institute of Mathematical Sciences, New York University. 
The authors are very grateful to Peter D. Lax for helpful remarks 
on a first draft of this paper. 
\endproclaim

%----------------------------------------------------------------
\heading{References} 
\endheading 

\ref \no\refBouchutBourdariasPerthame \by F. Bouchut, C. Bourdarias,
and B. Perthame
\paper A MUSCL methods satisfying 
all the numerical entropy inequalities
\jour Math. Comp. 
\vol\rm 65 \yr 1996 \pages 1439-1461
\endref

\ref   \no\refBrenierOsher  \by Y. Brenier and S.J. Osher
\paper The one-sided Lipschitz condition for convex scalar conservation laws   
\jour  SIAM J. Numer. Anal. \vol\rm 25 \yr 1988 \pages 8--23 
\endref

\ref   \no\refConwaySmoller \by E. Conway and J.A. Smoller 
\paper Global solutions of the Cauchy problem for quasilinear first order
 equations in several space variables 
\jour Comm. Pure Appl. Math. \vol\rm 19 \yr 1966 \pages 95--105
\endref

\ref \no\refCoquelLeFloch \by F. Coquel and P.G. LeFloch
\paper An entropy satisfying MUSCL scheme 
for systems of conservation laws 
\jour  Numer. Math. \vol\rm  74 \yr 1996 \pages   1--33    
\endref

\ref \no\refColella\by P. Colella 
\paper A direct Eulerian MUSCL scheme for gas dynamics
\jour SIAM J. Sci. Stat. Comput. \vol\rm 6 \yr 1985 \pages  104--117    
\endref

\ref \no\refCrandallMajda  \by M.G. Crandall and A. Majda 
\paper Monotone difference approximations for scalar conservation laws 
\jour Math. of Comp. \vol\rm 34 \yr 1980 \pages 1--21
\endref

\ref \no\refDafermos  \by C.M. Dafermos 
\paper Generalized characteristics in hyperbolic conservation laws: a
study of the structure and the asymptotic behavior of solutions, 
\jour in ``Nonlinear Analysis and Mechanics: Heriot-Watt symposium'',
ed. R.J. Knops, Pitman, London \vol\rm 1 \yr 1977 \pages 1--58   
\endref

\ref \no\refEngquistOsher \by B. Engquist and S.J. Osher 
\paper One-sided difference approximations for nonlinear conservation
laws 
\jour Math. of Comp. \vol\rm 35 \yr 1981 \pages 321--351
\endref

\ref \no\refFilippov  \by A.F. Filippov 
\paper Differential equations with discontinuous right-hand side 
\jour Mat. USSR Sbornik \vol\rm 51 \yr 1960 \pages 99--128 
\endref

\ref\no\refGlimmLax  \by J. Glimm and P.D. Lax 
\paper Decay of solutions of nonlinear hyperbolic conservation laws 
\jour Mem. Amer. Math. Soc \vol\rm 101 \yr 1970
\endref

\ref \no\refGodlevskiRaviart \by E. Godlevski and P.-A. Raviart 
\paper Hyperbolic Systems of Conservation Laws 
\jour Collection Ellipse, Paris \vol\rm  \yr 1993
\endref

\ref \no\refGoodmanLeVeque  \by J. Goodman and R. LeVeque
\paper A geometric approach to high resolution TVD schemes 
\jour SIAM J. Numer. Anal. \vol\rm 25 \yr 1988 \pages 268--284
\endref

\ref \no\refHartenone  \by A. Harten 
\paper On a class of high order resolution total-variation stable
finite difference schemes 
\jour SIAM J. Numer. Anal. \vol\rm 21 \yr 1974 \page 1--23 
\endref

\ref \no\refHartentwo  \by A. Harten 
\paper High resolution schemes for hyperbolic conservation laws 
\jour J. Comput. Phys. \vol\rm 49 \yr 1983 \page 357--393 
\endref

\ref \no\refHartenEngquistOsherChakravarthy
\by A. Harten, B. Engquist, S. Osher and S. Chakravarthy
\paper Uniformly high order accurate essentially non-oscillatory schemes
\jour J. Comput. Phys. \vol\rm 71 \yr 1987 \page 231--303 
\endref

\ref \no\refHartenHymanLax \by A. Harten, J.M. Hyman and P.D. Lax
\paper On finite-difference approximations and entropy conditions for
shocks  
\jour Comm. Pure Appl. Math.\vol\rm 29 \yr 1976 \pages 297--322 
\endref 

\ref \no\refJiangShu \by  G. Jiang and C.-W. Shu 
\paper On cell entropy inequality for 
discontinuous Galerkin methods 
\jour Math. Comp. 
\vol \rm 62 \yr 1994 \pages 531-538
\endref 

\ref \no\refKeyfitz  \by B.L. Keyfitz
\paper Solutions with shocks, an example of an $L^1$ contractive semigroup 
\jour Comm Pure Appl. Math. \vol\rm 24 \yr 1971 \pages 125--132
\endref

\ref \no\refKruzkov \by S.N. Kru\u zkov 
\paper First order quasilinear equations in several independent variables 
\jour Mat. USSR Sbornik \vol\rm 10 \yr 1970 \pages 217--243
\endref

\ref \no\refLaxone  \by P.D. Lax
\paper Hyperbolic systems of conservation laws II
\jour Comm. Pure Appl. Math. \vol\rm 10 \yr 1957 \pages 537--566
\endref

\ref \no\refLaxtwo  \by P.D. Lax
\book Hyperbolic Systems of Conservation Laws and the Mathematical
Theory of Shock Waves 
\publ SIAM \publaddr Philadelphia \yr 1973 
\endref

\ref \no\refLaxWendroff \by P.D. Lax and B. Wendroff
\paper Systems of conservation laws
\jour  Comm. Pure Appl. Math. \vol\rm 13 \yr 1960 \pages 217--237 
\endref

\ref \no\refvanLeerone \by B. van Leer
\paper Towards the ultimate conservative difference schemes, II, 
Monotonicity and conservation combined in a second order scheme
\jour J. Comp. Phys. \vol\rm 14 \yr 1974 \pages 361--370
\endref

\ref \no\refvanLeertwo \by B. van Leer
\paper Towards the ultimate conservative difference schemes, V, 
A second order sequel to Godunov's method
\jour J. Comp. Phys. \vol\rm 43 \yr 1981 \pages 357--372
\endref

\ref  \no\refLeFlochLiu  \by P.G. LeFloch and J.G. Liu
\paper   Discrete entropy and monotonicity criterion 
for hyperbolic conservation laws
\jour C.R. Acad. Sc. Paris \vol\rm 319 \yr 1994 \pages 881--886
\endref

\ref  \no\refLeroux  \by A.Y. Leroux
\paper Approximation de quelques probl\`emes hyperboliques 
nonlin\'eraires 
\jour Th\`ese d'\'etat, Universit\'e de Rennes \vol\rm \yr  1979 \pages 
\endref

\ref \no\refLionsSougadinisone \by P.-L. Lions and P. Sougadinis
\paper Convergence of MUSCL type schemes 
for scalar conservation laws  
\jour  C.R. Acad. Sc. Paris, S\'erie I, \vol\rm 311 \yr 1990 
\pages 259--264 
\endref

\ref \no\refLionsSougadinistwo \by P.-L. Lions and P. Sougadinis
\paper Convergence of MUSCL and filtered schemes 
for scalar conservation laws and Hamilton-Jacobi equations 
\jour  Numer. Math. \vol\rm 69 \yr 1995
\pages 441--470 
\endref

\ref \no\refNessyahuTadmor \by H. Nessyahu and E. Tadmor 
\paper The convergence rate of approximate solutions for nonlinear
scalar conservation laws
\jour SIAM J. Numer. Anal. \vol\rm 29 \pages 1505--1519 \yr 1992 
\endref

\ref \no\refNTT 
\by H. Nessyahu, E. Tadmor, and T. Tamir 
\paper The convergence rate of Godunov type schemes
\jour SIAM J. Numer. Anal. \vol\rm 31 \pages 1-16 \yr 1994
\endref 

\ref \no\refOsherone  \by S.J. Osher
\paper Riemann solvers, the entropy condition and difference
approximations 
\jour SIAM J. Numer. Anal. \vol\rm 21 \yr 1984 \pages 217--235
\endref

\ref \no\refOshertwo  \by S.J. Osher
\paper Convergence of generalized MUSCL schemes
\jour SIAM J. Numer. Anal. \vol\rm 22 \yr 1985 \pages 947--961 
\endref

\ref \no\refOsherTadmor \by S.J. Osher and E. Tadmor
\paper On the convergence of difference approximations 
to scalar conservation laws
\jour Math. of Comp. \vol\rm 50 \yr 1988 \pages 19--51 
\endref

\ref \no\refShu  \by C.W. Shu
\paper Total variation diminishing time discretizations
\jour SIAM J. Sci. Stat. Comp. \vol\rm 9 \yr 1988 \pages 1073--1084
\endref

\ref \no\refTadmorone   \by E. Tadmor
\paper Numerical viscosity and the entropy condition for conservative
difference schemes
\yr 1984 \jour Math. of Comp. \vol\rm 43 \pages 369--382 
\endref

\ref \no\refTadmortwo  \by E. Tadmor
\paper Convenient total variation diminishing conditions for nonlinear
difference schemes
\yr 1988 \jour SIAM J. Numer. Anal. \vol\rm 25 \pages 1002--1014 
\endref

\ref \no\refVolpert \by A.I. Volpert
\paper The space BV and quasilinear equations
\jour Mat. USSR Sbornik \vol\rm 73 \yr 1967 \pages 255--302
\endref

\ref \no \refYangone \by H. Yang
\paper On wavewise entropy inequalities for high-resolution 
schemes I~: the semidiscrete case 
\jour Math. of Comput. \vol\rm 213 
\yr 1996  \pages 45--67 \endref

\ref \no \refYangtwo \by H. Yang
\paper On wavewise entropy inequalities for high-resolution 
schemes II~: fully-discrete MUSCL schemes with exact evolution in 
small time 
\jour SIAM J. Numer. Anal., to appear 
\yr  \endref

\ 

%--------------------------------------------------------------
\rm 

\subheading{Appendix: Monotonicity Property}

In this appendix we brieffly discuss the monotonicity property together
with more basic properties of entropy solutions to 
conservation laws, which go back to Kruzkov \cite{\refKruzkov}
and Volpert \cite{\refVolpert}. 
In the paper by Keyfitz, \cite{\refKeyfitz}, 
somewhat simpler proofs are available for piecewise Lipschitz continuous
solutions. We are interested in the local versions of the properties,
i.e.~formulated in domains limited by characteristic curves. 
To cope with discontinuous solutions, we use 
the concept of generalized characteristic curves introduced for 
ordinary differential equations by Filippov \cite{\refFilippov} 
and developed in the context of conservation laws
by Dafermos; see e.g.~\cite{\refDafermos} and the references therein. 
We recall that, through any point $(x_0,t_0)$, 
there exists a funnel of forward
and backward generalized characteristic curves, which 
fill up a domain $\big\{\xi^m(t) \le x \le \xi^M(t)\big\}$. 
Here $\xi^m$ (respectively $\xi^M$) is called the minimal (resp. maximal) 
characteristic curve originating at $(x_0,t_0)$.
It is known \cite{\refFilippov, \refDafermos} that a
characteristic, say $\xi$, is Lipschitz continuous and 
for almost every time $t>0$ satisfies
$$
\frac{d \xi}{dt} (t) = \cases f'(u(\xi(t),t)) 
& \text{if } u_-(t) = u_+(t) = u(\xi(t),t),\\
        \frac{f(u_+)-f(u_-)}{u_+-u_-} &\text{if } u_-\ne u_+, 
\endcases
\tag A.1
$$
where $u_{\pm}=u\bigl(\xi(t)\pm,t\bigr)$.
For our purposes, $f$ is strictly convex and 
there is a unique forward  characteristic issued from  
$(x_0),t_0$ and there is no need to distinguish between the minimal
characteristic and the maximal one,
with the exception of those points where $t_0=0$ and $u_0$ has an
increasing jump at $x_0$; cf.~\cite{\refDafermos}.

Solutions $u$ to (1.1) are Lipschitz continuous in
time with values in $L^1$ and,  for all times $t$,
$u(t)$ has bounded total variation in $x$.

The following properties follows from \cite{\refVolpert, \refKruzkov} 
and the technique of generalized characteristic in 
\cite{\refFilippov, \refDafermos}.

\proclaim{Proposition A.1} 
Let $u$ be the entropy solution to $(1.1)$-$(1.2)$. 
Given $x_1$ and $x_2$ with $x_1 < x_2$,
consider the maximal forward characteristic $\xi_1^u(t)$ 
issued from $(0,x_1)$ and
the minimal forward characteristic $\xi_2^u(t)$ from $(0,x_2)$. 
For all times $t\ge0$, $u$  satisfies 
\roster
\item the local {\rm maximum principle}
for all $t \le s$ and $y \in (\xi_1^u(s), \xi_2^u(s))$:
$$
\inf_{\xi_1^u(t)< x < \xi_2^u(t)}  \big\{ u(x,t) \big\}
\le u(s,y) 
\le \sup_{\xi_1^u(t)< x < \xi_2^u(t)}  \big\{ u(x,t) \big\}, 
\tag A.2
$$
\item the local {\rm $L^1$ contraction property}: 
$$
\frac{d}{dt}\int_{\max(\xi_1^u(t),\xi_1^v(t))}^{\min(\xi_2^u(t),\xi_2^v(t))} 
\bigl| u(x,t) -v(x,t) \bigr| \,dx \le 0,
\tag A.3
$$
\item the local {\rm order preserving property}: 
$$
\aligned
&\text{ if }  \, \, u_0(x) \le v_0(x)  \,  \,\text{ for  all } \ \, \, \,
            x \in \bigl(\max(\xi_1^u(0),\xi_1^v(0)), \min(\xi_2^u(0),\xi_2^v(0))\bigr),\\ 
&\text{ then } \, \, u(x,t) \le v(x,t)  \, \,\text{ for all }  \, \, \,
             x \in \bigl(\max(\xi_1^u(t),\xi_1^v(t)) < x < \min(\xi_2^u(t),\xi_2^v(t))\bigr), 
\endaligned
\tag A.4
$$
\item and the local TVD {\rm property\/}: 
$$
\frac{d}{dt}\operatorname{TV}\limits_{ \xi_1^u(t)}^{ \xi_2^u(t)} u(t) \le 0.
\tag A.5
$$
\endroster
where in $(1)$ and $(3)$ the function $v$ denotes the solution to $(1.1)$-$(1.2)$ 
with $u_0$ replaced by a function $v_0$ $\in$ $BV(\RR)$, and the curves 
$\xi_1^v(t)$ and $\xi_2^v(t)$ are defined in the obvious way. $\hfill \square$
\endproclaim

In (A.5), $\operatorname{TV}\limits_{a}^{b} (w)$ denotes the
total variation of a function $w: (a,b) \to \RR$. The derivatives in
(A.3) and (A.5) are to be understood in the distributional sense.
At least with $x_1= \xi_1^u=-\infty$ and $x_2= \xi_2^u =\infty$, 
(A.5) is a direct consequence of (A.3) and the
invariance by translation of the solution-operator for (1.1).

The following proposition concerns the {\it monotonicity property\/},
which is a refinement to the statement that 
the solution-operator is {\it monotonicity preserving\/}, i.e. satisfies     
$$
\text{ if } \, \,  u_0 \, \,  \text{ is monotone, then } \, \, 
u(t) \, \,  \text{ is  monotone for all } t \ge 0.
\tag A.6
$$
Namely (A.6) follows from (A.3) by taking 
$x_1= \xi_1^u=-\infty$ and $x_2= \xi_2^u =\infty$
and $v_0(x) = u_0(x+y)$ for positive or negative values of $y$. 

It is convenient now to assume that (1.4) is satisfied. 

\proclaim{Proposition A.2} 
(Monotonicity property) Suppose that the initial
condition \, $u_0$ \, has a locally finite number of local extrema.
There exist (Lipschitz  continuous) generalized characteristic curves 
$\varphi_q: [\,0,\infty) \to \RR$ 
--the index $q$ describing a subset $E(u_0)$ of consecutive integers-- such that 
$$
\varphi_q \le \varphi_{q+1},
\tag A.7
$$
$$
\text{ there is only a finite number of such curves in each compact set}
\tag A.8 
$$ 
and, for all $t\ge0$ and all relevant values of $p$, 
$$
\aligned
&   u(t) \,\, \text{ is non-decreasing for } x \in
     (\varphi_{2p}(t),\, \varphi_{2p+1}(t)),\\
&   u(t) \,\, \text{ is non-increasing for } x \in 
    (\varphi_{2p-1}(t),\, \varphi_{2p}(t)),
\endaligned
\tag A.9
$$
and, as long as $\varphi_{2p}(t)\ne \varphi_{2p+1}(t)$,
$$
\aligned
&   u(\varphi_{2p}(t)+,t)  \, \, \text{ is non-decreasing, } \\
&   u(\varphi_{2p+1}(t)-,t)  \, \, \text { is non-increasing.}
\endaligned
\tag A.10
$$
$\hfill \square$
\endproclaim

The paths $\varphi_{2p}$ and $\varphi_{2p+1}$ are called 
path of local maximum and path of local minimum for $u$, respectively.
By convention, 
$\varphi_{q} \equiv -\infty$ for $q>\max E(u_0)$ and $\varphi_{q} \equiv +\infty$ 
for $q<\min E(u_0)$. 

If two initially distinct paths cross at a 
later time, then from that time they will coincide thanks to the uniqueness 
property for forward characteristics. Note also that such paths 
need no longer be paths of local extremum in a strict sense, 
but arbitrary characteristics , even though (A.9) would still hold. 
On the other hand, when simultaneously $u_0$
has a decreasing jump at a point $x_0$, that 
$u_0(x_0-)$ is a local maximum, and $u(x_0+)$ is a local minimum, 
then two equal paths originate from $(0, x_0)$, one being 
a path of minimum and the other a path of maximum. 
When simultaneously $u_0$ has an increasing jump at $x_0$, 
that $u_0(x_0-)$ is a local minimum, and that $u_0(x_0+)$ a local maximum, 
then two distinct paths originate from $(0, x_0)$. 

Proposition A.2 is a classical matter in the literature
although no specific reference seems available. Cf.~however 
Harten \cite{\refHartenone} and Tadmor \cite{\refTadmorone} where 
the ideas are developed. 

\proclaim{Proof of Proposition A.2} \rm 
First of all, the points $\varphi_q (0)$ and the set $E(u_0)$ are defined
from the initial condition $u_0$ in an obvious way so that the conditions
(A.7)--(A.10) hold true at time $t=0$. 
Let us define $\varphi_{2p}(t)$ to be
the maximal forward characteristic issued from $\varphi_{2p}(0)$.
Similarly, let $\varphi_{2p+1}(t)$ be the minimal forward
characteristics issued from $\varphi_{2p+1}(0)$. 
Indeed one need to distinguish between minimal and maximal characteristics 
only in the case of an initially increasing jump.  
It may happen that both a path of minimum
{\it and\/} a path of maximum may originate from such a point of
increasing jump. 

Property (A.7) is an immediate consequence of the uniqueness property of 
the forward characteristic. 
Condition (A.8) follows from the property of propagation with finite speed 
satisfied by solutions to (1.1)  and the fact that the initial data has a locally finite
number of local extremum. 
Indeed (A.1) yields a uniform bound for the slopes of the characteristics.

In order to establish (A.9) and (A.10), 
we first suppose that $u_0$ does not admit increasing jumps. 
Consider an interval of the form $(\varphi_{2p}(t), \varphi_{2p+1}(t))$ 
for those values of $t$ when this interval is not empty. 
Note first that, taking $\xi_1^u=\varphi_{2p}$ and
$\xi_2^u=\varphi_{2p+1}$, the local maximum principle (A.2) implies 
in particular that 
$$
\aligned
u\bigl(\varphi_{2p}(t)+, t\bigr) \ge u_0\bigl(\varphi_{2p}(0)+\bigr),\\
u\bigl(\varphi_{2p+1}(t)-,t\bigr) \le u_0\bigl(\varphi_{2p+1}(0)-\bigr).
\endaligned
\tag A.11
$$
Let $w$ be the entropy solution to (1.1) associated with the initial condition
$$
w(x,0) = w_0(x) \equiv
\aligned
&   u_0\bigl(\varphi_{2p}(0)+\bigr) \qquad \text{if } x < \varphi_{2p}(0), \\
&   u_0(x)    \qquad \text{if } \varphi_{2p}(0) < x < \varphi_{2p+1}(0), \\
&   u_0\bigl(\varphi_{2p+1}(0)-\bigr) \qquad  \text{if } x > \varphi_{2p+1}(0).
\endaligned
\tag A.12
$$
The data $w_0$ is non-decreasing and, in view of (A.6),
$$
\text{the solution }  w \text{ is non-decreasing for all times.}
\tag A.13
$$ 
Let $\psi_{2p}$ and $\psi_{2p+1}$ 
be the forward characteristics associated with $w$ and issued  
$\varphi_{2p}(0)$ and $\varphi_{2p+1}(0)$ at time $t=0$, respectively. 
Observe that the maximum forward and the 
minimum forward curves coincide since by construction $w_0$ is continuous 
at $\varphi_{2p}(0)$ and $\varphi_{2p+1}(0)$.
Note in passing that the function $w$ satisfies: 
$$
w(x,t) = 
\cases u_0\bigl(\varphi_{2p}(0)+\bigr) & \text{if } x < \psi_{2p}(t),\\
u_0\bigl(\varphi_{2p+1}(0)-\bigr) & \text{if } x > \psi_{2p+1}(t).
\endcases
\tag A.14
$$

Using (A.11) and (A.14) and the fact that $f'(.)$ is increasing, one gets 
$$
\frac{d \varphi_{2p+1}}{dt} (t) \le f'(u(\varphi_{2p+1}(t)-,t)) 
\le f'(u_0(\varphi_{2p+1}(0)-)) 
\le \frac{d \psi_{2p+1}}{dt} (t), 
$$
which implies 
$$
\varphi_{2p+1} (t) \le \psi_{2p+1} (t).
$$
Similarly
$$
\psi_{2p} (t) \le \varphi_{2p}(t).
$$
Using the $L^1$ contraction principle (A.3), it follows that
$$
u = w  \, \qquad \text{ for } \,  { \varphi_{2p}(t) < x < \varphi_{2p+1}(t) }, 
$$
and, in view of (A.13), the function $u$ is non-decreasing and (A.9) holds. 
Using (A.9) and the local maximum principle (A.2) finally provides (A.10). 
The proof is complete in the case of an interval 
of the form $(\varphi_{2p}(t), \varphi_{2p+1}(t))$.
An interval  $(\varphi_{2p-1}(t), \varphi_{2p}(t))$ can be treated in a similar fashion. 

It remains to consider increasing jumps in $u_0$. That situation
can be treated by using the following property.
Suppose $u_0$ has an increasing jump at a point $x_0$ and let 
$\varphi^m(t)$ and $\varphi^M(t)$ be the minimal and maximal forward
curves from $x_0$. It is known that at least for small times
the function $u(t)$ coincides with the rarefaction wave connecting the 
values $u_0(x_0\pm)$ in the interval $\bigl(\varphi^m(t),\varphi^M(t)\bigr)$.
 
This completes the proof of Proposition A.2.        $\hfill \square$
\endproclaim

\enddocument